\documentclass{jotart}
\usepackage{amsmath}
\usepackage{amsxtra}
\usepackage{amsfonts}
\usepackage{amssymb}
\theoremstyle{proclaim}
\newtheorem{theorem}{Theorem}[section]
\newtheorem{corollary}[theorem]{Corollary}
\newtheorem{lemma}[theorem]{Lemma}
\newtheorem{proposition}[theorem]{Proposition}
\theoremstyle{statement}
\newtheorem{definition}[theorem]{Definition}
\newtheorem{remark}[theorem]{Remark}
\newtheorem{example}[theorem]{Example}

\numberwithin{equation}{section}

\def\hooklongrightarrow{\lhook\joinrel\longrightarrow}
\newcommand{\ip}[2]{\left\langle\, #1 \mid #2 \,\right\rangle}
\newcommand{\bip}[2]{\left\langle\, #1 \bigm| #2 \,\right\rangle}
\newcommand{\Bip}[2]{\left\langle\, #1 \Bigm| #2 \,\right\rangle}
\newcommand{\cj}[1]{\overline{#1}}
\newcommand{\bz}{\mathbb{Z}}
\newcommand{\br}{\mathbb{R}}
\newcommand{\bc}{\mathbb{C}}
\newcommand{\bt}{\mathbb{T}}

\newbox\frogdown
\newlength\frogdrop
\def\hookdownarrow{\setlength{\unitlength}{0.4pt}\setbox\frogdown=\hbox to 0pt{\hss $\displaystyle \downarrow $\hss }\setlength{\frogdrop}{2.5\ht\frogdown}\raisebox{0pt}[5\unitlength][\frogdrop]
{\begin{picture}(0,5)(10,0)
\put(5,0){\oval(10,10)[t]}
\end{picture}\lower\ht\frogdown\box\frogdown}}

\begin{document}
\issueinfo{00}{0}{0000} 
\commby{S. Stratila}
\pagespan{101}{144}
\date{November dd, 2005}
\title[Martingales and covariant systems]{Martingales, endomorphisms, and covariant systems of operators in Hilbert
space}
\dedicatory{Dedicated to the memory of J.L. Doob}
\author[Dutkay {\protect \and} Jorgensen]{Dorin Ervin Dutkay {\protect
\and} Palle E.T. Jorgensen}
\address{\uppercase{Dorin Ervin Dutkay}, Department of Mathematics, Hill Center-Busch Campus, Rutgers, The State University of New Jersey, Piscataway, NJ 08854-8019, U.S.A.}
\email{ddutkay@math.rutgers.edu}
\address{\uppercase{Palle E.T. Jorgensen}, Department of Mathematics, The University of Iowa, Iowa City, IA 52242-1419, U.S.A.}
\email{jorgen@math.uiowa.edu}
\begin{abstract} 
In the theory of wavelets, in the study of subshifts, in the
analysis of Julia sets of rational maps of a complex variable,
and, more generally, in the study of dynamical systems, we are
faced with the problem of building a unitary operator from a
mapping $r$ in a compact metric space $X$. The space $X$ may be a
torus, or the state space of subshift dynamical systems, or a
Julia set.
\par
 While our motivation derives from some wavelet problems, we have in mind other applications as well; and the issues involving
covariant operator systems may be of independent interest.
\end{abstract}
\maketitle

\tableofcontents
\section{\label{Intro}\uppercase{Introduction}}

      In this paper, we aim at combining and using ideas from one area of
mathematics (operator theory and traditional analysis) in a different area
(martingale theory from probability). We have in mind applications to both
wavelets and symbolic dynamics. So our paper is interdisciplinary:
results in one area often benefit the other. In fact, the benefits go both
ways.

\par
 Our construction is based on a closer examination of an eigenvalue problem
for a transition operator, also called a Perron-Frobenius-Ruelle
operator.
\par
Under suitable conditions on the given filter functions, our
construction takes place in the Hilbert space $L^2(\br^d)$. In a
variety of examples, for example for frequency localized wavelets,
more general filter functions are called for. This then entails
basis constructions in Hilbert spaces of $L^2$-martingales. These
martingale Hilbert spaces consist of $L^2$ functions on certain
projective limit spaces $X_\infty$ built on a given mapping $r:
X\rightarrow X$ which is onto, and finite-to-one. We study
function theory on $X_\infty$ in a suitable general framework, as
suggested by applications; and we develop our theory in the
context of Hilbert space and operator theory.

      We hope that these perhaps unexpected links between more traditional
and narrowly defined fields will inspire further research. Since we wish to
reach several audiences, we have included here a few more details than is
perhaps standard in more specialized papers. The general question we address
already has a number of incarnations in the literature, but they have so far
not been unified. Here are two such examples which capture the essence of
our focus. (a)~Extension of non-invertible endomorphisms in one space $X$ to
automorphisms in a bigger space naturally containing $X$. (b)~Some
non-invertible operator $S$ (contractive or isometric) in a fixed Hilbert
space $\mathcal{H}$ is given. It is assumed that $S$ is contractive and that it satisfies
a certain covariance condition specified by a system of operators in $\mathcal{H}$. The
question is then to extend $S$ to a unitary operator $U$ in a bigger Hilbert
space which naturally contains $\mathcal{H}$, such that $U$ satisfies a covariance
condition arising by dilation from the initially given system on $\mathcal{H}$.
\par
The dilation idea in operator theory is fundamental; i.e., the
idea of extending (or dilating) an operator system on a fixed Hilbert space
$H_0$ to a bigger ambient Hilbert space $H$ in such a way as to get
orthogonality relations in the dilated space $H$; see for example \cite{PaSc72} and
Remark \ref{rem3_2} below.  In an operator algebraic framework such an
extension is of course encoded by Stinespring's theorem \cite{Sti55}.  Our
present setting is motivated by this, but goes beyond it in a number of
ways, as we show in Sections \ref{ProLim}--\ref{Ite} below. 
\par
          Our basic viewpoint may be understood from the example of
wavelets: A crucial strength of wavelet bases is their algorithmic and
computational features. What this means in terms of the two Hilbert spaces
are three things: First we must have a concrete function representation of
the dilated space $H$; and secondly we aim for recursive and matrix based
algorithms, much like the familiar case of Gram-Schmidt algorithms which
lets us compute orthonormal bases, or frames (see e.g., \cite{BJMP05}) in the
dilated space $H$. Thirdly, we reverse the traditional point of view. Hence,
the dilation idea is turned around: Starting with $H$, we wish to select a
subspace $H_0$ which is computationally much more feasible. This idea is
motivated by image processing where such a selected subspace $H_0$ corresponds
to a chosen resolution, and where "resolution" is to be understood in the
sense of optics; see e.g., \cite{JMR01} and \cite{Jor06}. The selection of subspace
$H_0$ is made in such a way as to yield recursive algorithms to be used in
computation of orthonormal bases, or frames in $H$, but starting with data
from $H_0$.
\par
      Examples of (a) occur in thermodynamics, such as it is presented in
its rigorous form by David Ruelle in \cite{Rue89} and \cite{Rue04}. Both (a) and (b)
are present in the approach to wavelets that goes under the name
multiresolution analysis (MRA) \cite{Dau92}. In this case, we can take $X$ to be
$\mathbb{R}/\mathbb{Z}$, or equivalently the circle, or the unit-interval $[0,1)$, and the
extension of $X$ can be taken to be the real line $\mathbb{R}$ (see \cite{Dau92}), or it may
be a suitable solenoid over $X$; see, e.g., \cite{BreJo91} and \cite{Bre96b}. In this case,
the endomorphism in $X$ is multiplication by $2$ modulo the integers $\mathbb{Z}$, and the
extension to $\mathbb{R}$ is simply $x \to 2x$. The more traditional settings for (b)
are scattering theory \cite{LaPh76} or the theory of extensions, or unitary
dilations of operators in Hilbert space, as presented for example in
\cite{JoMu80}, \cite{BMP00}, and in the references given there.

Specifically, we study the problem of inducing operators on
Hilbert space from non-invertible transformations on compact
metric spaces. The operators, or representations must satisfy
relations which mirror properties of the given point
transformations.
\par
While our setup allows a rather general formulation in the context
of $C^*$-algebras, we will emphasize the case of induction from an
abelian $C^*$-algebra. Hence, we will stress the special case when
$X$ is a given compact metric space, and $r:X\rightarrow X$ is a
finite-to-one mapping of $X$ onto $X$. Several of our results are
in the measurable category; and in particular we are not assuming
continuity of $r$, or any contractivity properties.
\subsection{\label{Wav}Wavelets}
\par Our results will apply to wavelets. In the theory
of multiresolution wavelets, the problem is to construct a special
basis in the Hilbert space $L^2(\br^d)$ from a set of numbers
$a_n$, $n\in\bz^d$.
\par
The starting point is the scaling identity
\begin{equation}\label{eqintsc1}
\varphi(t)=N^{1/2}\sum_{n\in\bz^d}a_n\varphi(At-n),\quad(t\in\br^d),
\end{equation}
where $A$ is a $d$ by $d$ matrix over $\bz$, with eigenvalues
$|\lambda|> 1$, and $N=|\mbox{det} A|$, and where $\varphi$ is a
function in $L^2(\br^d)$.
\par
The first problem is to determine when (\ref{eqintsc1}) has a
solution in $L^2(\br^d)$, and to establish how these solutions
({\it scaling functions}) depend on the coefficients $a_n$.
\par
When the Fourier transform is applied, we get the equivalent
formulation,
\begin{equation}\label{eqintsc2}
\hat\varphi(x)=N^{-1/2}m_0({A^{tr}}^{-1}x)\hat\varphi({A^{tr}}^{-1}x),
\end{equation}
where $\hat\varphi$ denotes the Fourier transform,
\begin{equation*}\hat\varphi(x)=\int_{\br^d}\mathrm{e}^{-i2\pi x\cdot t}\varphi(t)\,\mathrm{d}t\end{equation*}
and where now $m_0$ is a function on the torus
\begin{equation*}\bt^d=\{z=(z_1,...,z_d)\in\bc^d\,
|\, |z_j| =1, 1\leq j\leq d\}=\br^d/\bz^d\end{equation*}, i.e.,
\begin{equation*}m_0(z)=\sum_{n\in\bz^d}a_nz^n=\sum_{n\in\bz^d}a_ne^{-i2\pi
n\cdot x}.\end{equation*} The duality between the compact group $\bt^d$ and the
lattice $\bz^d$ is given by
\begin{equation*}\ip{z}{n}=z^n=z_1^{n_1}...z_d^{n_d},\quad(z=(z_1,...,z_d),n=(n_1,...,n_d)).\end{equation*}
\par In this case, matrix multiplication $x\mapsto A x$ on $\br^d$ passes to the quotient $\br^d/\bz^d$, and we get an $N$-to-one
mapping $x\mapsto Ax\mod\bz^d$, which we denote by $r_A$.
\par
The function $m_0$ is called a low pass filter, and it is chosen
such that the operator $S=S_{m_0}$ given by
\begin{equation*}(Sf)(z)=m_0(z)f(Az)\end{equation*}
is an isometry on $H_0=L^2(\bt^d,\mbox{Haar measure})$. Moreover,
$L^\infty(\bt^d)$ acts as multiplication operators on $H_0$. If
$g\in L^\infty(\bt)$
\begin{equation*}(M(g)f)(z)=g(z)f(z)\end{equation*}
and
\begin{equation}\label{eqintcov1}
SM(g)=M(g(A\cdot))S
\end{equation}
A main problem is the extension of this covariance relation
(\ref{eqintcov1}) to a bigger Hilbert space $H_0\rightarrow
H_{ext}$, $S\rightarrow S_{ext}$, such that $S_{ext}$ is unitary
in $H_{ext}$. We now sketch briefly this extension in some
concrete cases of interest.
\par
In Section \ref{ProLim}, we construct a sequence of measures
$\omega_0,\omega_1,...$ on $\bt^d$ such that
$L^2(\bt^d,\omega_0)\simeq H_0$, and such that there are natural
isometric embeddings
\begin{equation}\label{eqint3}L^2(\bt^d,\omega_n)\hookrightarrow
L^2(\bt^d,\omega_{n+1}),\quad f\mapsto f\circ r_A.\end{equation}
The limit in (\ref{eqint3}) defines a {\it martingale Hilbert
space} $\mathcal{H}$ in such a way that the norm of the
$L^2$-martingale $f$ is
\begin{equation*}\|f\|^2=\lim_{n\rightarrow\infty}\|P_nf\|^2_{L^2(\bt^d,\omega_n)}\end{equation*}
We also state a pointwise a.e.\ convergence result (Section \ref{Mar}).
\par
If $\Psi: L^2(\bt^d,\omega_n)\rightarrow L^2(\br^d)$ is defined by
\begin{equation*}\Psi: f_n\mapsto f_n(A^{-n}x)\hat\varphi(x),\end{equation*}
then $\Psi$ is an isometry of $L^2(\bt^d,\omega_n)$ into
$L^2(\br^d)$.
\par
Specifically
\begin{equation}\label{eqint4}
\int_{\bt^d}|f_n|^2\,\mathrm{d}\omega_n=\int_{\br^d}|f_n(A^{-n}x)\hat\varphi(x)|^2\,\mathrm{d}x
\end{equation}
As a result we have induced a system
\begin{equation*}(r_A,\bt^d)\rightarrow (S_{m_0},L^2(\bt^d))\rightarrow
(U_A,L^2(\br^d)).\end{equation*} where
\begin{equation}\label{eqintdil}
(U_A\xi)(x)=N^{1/2}f(Ax),\quad(f\in L^2(\br^d)) \end{equation}
$U_A$ unitary; the system is determined by the given filter
function $m_0$. It can be checked (see details in Section \ref{Mar}) that
$\Psi$ is an isometry, and that
\begin{equation*}U_AM(g)=M(g(A\cdot))U_A\end{equation*}
holds on $L^2(\br^d)$. Moreover $\Psi$ maps onto $L^2(\br^d)$ if
the function $m_0$ doesn't vanish on a subset of positive measure.
\par
In the case of wavelets, we ask for a wavelet basis in
$L^2(\br^d)$ which is consistent with a suitable resolution
subspace in $L^2(\br^d)$. Whether the basis is orthonormal, or
just a Parseval frame, it may be constructed from a system of
subband filters $m_i$, say with $N$ frequency bands. These filters
$m_i$ may be realized as functions on $X=\bt^d=\br^d/\bz^d$, the
$d$-torus. Typically the scaling operation is specified by a given
expansive integral $d$ by $d$ matrix $A$.
\par
 Let $N:= |\mbox{det} A|$. Pass $A$ to the quotient $X$ = $\br^d/\bz^d$, and we get a mapping
$r$ of $X$ onto $X$ such that $\# r^{-1}(x) = N$ for all $x$ in
$X$, and the $N$ branches of the inverse are strictly contractive
in $X=\br^d/\bz^d$ if the eigenvalues of $A$ satisfy
$|\lambda|>1$.
\par
The subband filters $m_i$ are defined in terms of this map, $r_A$,
and the problem is now to realize the wavelet data in the Hilbert
space $L^2(\br^d)$ in such a way that $r=r_A :X\rightarrow X$
induces the unitary scaling operator $f\mapsto N^{1/2} f(A x)$ in
$L^2(\br^d)$, see (\ref{eqintdil}).
\subsection{\label{Jul}Examples (Julia sets, subshifts)}
\par
In this paper we will show that this extension from spaces $X$,
with a finite-to-one mapping $r:X\rightarrow X$, to operator
systems may be done quite generally, to apply to the case when $X$
is a Julia set for a fixed rational function of a complex
variable, i.e., $r(z)=p_1(z)/p_2(z)$, with $p_1,p_2$ polynomials,
$z\in\bc$ and $N=\max(\mbox{deg}\,p_1,\mbox{deg}\,p_2)$. Then
$r:X(r)\rightarrow X(r))$ is $N$-to-$1$ except at the singular
points of $r$. Here $X(r)$ denotes the Julia set of $r$. \par
 It
also applies to shift invariant spaces $X(A)$ when $A$ is a $0-1$
matrix, and
\begin{equation*}X(A)=\{(x_i)\in\prod_{\mathbb{N}}\{1,...,N\}\,|\,A(x_i,x_{i+1})=1\}\end{equation*}
and
\begin{equation*}r_A(x_1,x_2,...)=(x_2,x_3,...)\end{equation*}
is the familiar subshift. Note that $r_A:X(A)\rightarrow X(A)$ is
onto iff every column in $A$ contains at least one entry $1$.
\subsection{\label{SMar}Martingales}
\par
Part of the motivation for our paper derives from the papers by
Richard Gundy \cite{Gun00}, \cite{Gun04}, \cite{Gun99},
\cite{Gun66}. The second named author also acknowledges
enlightening discussions with R. Gundy. The fundamental idea in
these papers by Gundy et al is that multiresolutions should be
understood as martingales in the sense of Doob
\cite{Doob1},\cite{Doob2},\cite{Doob3} and Neveu \cite{Neveu}. And
moreover that this is a natural viewpoint.
\par
  One substantial advantage of this viewpoint is that we are then able to
handle the construction of wavelets from subband filters that are
only assumed measurable, i.e., filters that fail to satisfy the
regularity conditions that are traditionally imposed in wavelet
analysis.
\par
  A second advantage is that the martingale approach applies to a number of
wavelet-like constructions completely outside the traditional
scope of wavelet analysis in the Hilbert space $L^2(\br^d)$. But
more importantly, the martingale tools apply even when the
operation of scaling doesn't take place in $\br^d$ at all, but
rather in a compact Julia set from complex dynamics; or the
scaling operation may be one of the shift in the subshift dynamics
that is understood from that thermodynamical formalism of David
Ruelle \cite{Rue89}.
\subsection{\label{Gen}The general theory}
\par
 In each of the examples, we are faced with a given space $X$, and a
finite-to-one mapping $r:X\rightarrow X$. The space $X$ is
equipped with a suitable family of measures $\mu_h$, and the
$L^\infty$ functions on $X$ act by multiplication on the
corresponding $L^2$ spaces, $L^2(X,\mu_h)$. It is easy to see that
there are $L^2$ isometries which intertwine the multiplication
operators $M(g)$ and $M(g\circ r)$, as $g$ ranges over
$L^\infty(X)$. We have
\begin{equation}\label{eqintadd1}
\begin{array}{ccc}
L^2(X,\mu_h) & \smash{\overset{\displaystyle S}{\longrightarrow}} & L^2(X,\mu_h) \\
\hookdownarrow & & \hookdownarrow \\
\mathcal{H}_{ext} & \smash{\overset{\displaystyle U}{\longrightarrow}} & \mathcal{H}_{ext}
\end{array}
\end{equation}
where the vertical maps are given by inclusions. Specifically,
\begin{equation}\label{eqintadd2}
SM(g)=M(g\circ r)S,\mbox{ and } UM(g)U^{-1}=M(g\circ r)
\end{equation}
\par
But for spectral theoretic calculations, we need to have
representations of $M(g)$ and $M(g\circ r)$ unitarily equivalent.
That is true in traditional wavelet applications, but the unitary
operator $U$ in (\ref{eqintadd2}) is not acting on $L^2(X,\mu_h)$.
Rather, the unitary $U$ is acting by matrix scaling on a different
Hilbert space, namely\\ $L^2(\br^d,\mbox{Lebesgue measure})$,
\begin{equation*}U_Af(t)=|\mbox{ det }A|^{1/2}f(At),\quad(t\in\br^d,f\in
L^2(\br^d).\end{equation*}
\par
In the other applications, Julia set, and shift-spaces, we aim for
a similar construction. But in these other cases, it is not at all
clear what the Hilbert space corresponding to $L^2(\br^d)$, and
the corresponding unitary matrix scaling operator, should be.
\par
We provide two answers to this question, one at an abstract level,
and a second one which is a concrete function representation;
Sections \ref{Mul} and \ref{ProLim}.
\par
At the abstract level, we show that the construction may be
accomplished in Hilbert spaces which serve as unitary dilations of
the initial structure, see (\ref{eqintadd1}). \par In the
concrete, we show that the extended Hilbert spaces may be taken as
Hilbert spaces of $L^2$- martingales on $X$. In fact, we present
these as Hilbert spaces of $L^2$ functions built from a projective
limit
\begin{equation*}X\stackrel{r}{\leftarrow}X\stackrel{r}{\leftarrow}X....\leftarrow
X_\infty.\end{equation*}
 This is analogous to the distinction between an
abstract spectral theorem on the one hand, and a concrete spectral
representation, on the other. To know details about
multiplicities, and multiplicity functions (Section \ref{Mul}), we need
the latter.
\par
Our concrete version of the dilation Hilbert space
$\mathcal{H}_{ext}$ from (\ref{eqintadd1}) is then
\begin{equation*}\mathcal{H}_{ext}\simeq L^2(X_\infty,\hat\mu_h)\end{equation*}
for a suitable measure $\hat\mu_h$ on $X_\infty$.
\section{\label{Fun}\uppercase{Functions and measures on} $X$}
Consider
\begin{itemize}
\item $X$ a compact metric space, \item
$\mathfrak{B}=\mathfrak{B}(X)$ a Borel sigma-algebra of subsets of
$X$, \item $r:X\rightarrow X$ an onto, measurable map such that
$\# r^{-1}(x)<\infty$ for all $x\in X$, \item
$W:X\rightarrow[0,\infty)$, \item $\mu$ a positive Borel measure
on $X$.
\end{itemize}
\subsection{\label{Tra}Transformations
of functions and measures}
\begin{itemize}
\item Let $g\in L^\infty(X)$. Then
\begin{equation}\label{eqfm2_2_1}
M(g)f=gf
\end{equation}
is the multiplication operator on $L^\infty(X)$ or on
$L^2(X,\mu)$. \item Composition:
\begin{equation}\label{eqfm2_2_2}
S_0f=f\circ r,\mbox{ or }(S_0f)(x)=f(r(x)),\quad(x\in X).
\end{equation}
\item If $m_0\in L^\infty(X)$, we set
\begin{equation*}S_{m_0}=M(m_0)S_0,\end{equation*}
or equivalently
\begin{equation}\label{eqfm2_2_3}
(S_{m_0}f)(x)=m_0(x)f(r(x)),\quad(x\in X,f\in L^\infty(X)).
\end{equation}
\item $r^{-1}(E):=\{x\in X\,|\, r(x)\in E\}$ for
$E\in\mathfrak{B}(X)$.
\begin{equation*}\mu\circ r^{-1}(E)=\mu(r^{-1}(E)),\quad(E\in\mathfrak{B}(X)).\end{equation*}
\end{itemize}

\subsection{\label{ProMea}Properties of measures $\mu$ on $X$. Definitions}
\begin{enumerate}
\item {\it Invariance}:
\begin{equation}\label{eqfm4}
\mu\circ r^{-1}=\mu.
\end{equation}
\item {\it Strong invariance}:
\begin{equation}\label{eqfm5}
\int_Xf(x)\,\mathrm{d}\mu=\int_X\frac{1}{\#r^{-1}(x)}\sum_{r(y)=x}f(y)\,\mathrm{d}\mu,\quad(f\in
L^\infty(X)).
\end{equation}
\item $W:X\rightarrow[0,\infty)$,
\begin{equation}\label{eqfm6}
(R_Wf)(x)=\sum_{r(y)=x}W(y)f(y).
\end{equation}
If $m_0\in L^\infty(X,\mu)$ is complex valued, we use the notation
$R_{m_0}:=R_W$ where $W(x)=|m_0(x)|^2/\#r^{-1}(r(x))$.
\begin{enumerate}
\item A function $h:X\rightarrow[0,\infty)$ is said to be an
eigenfunction for $R_W$ if
\begin{equation}\label{eqfm7}
R_Wh=h
\end{equation}
\item A Borel measure $\nu$ on $X$ is said to be a
left-eigenfunction for $R_W$ if
\begin{equation}\label{eqfm8}
\nu R_W=\nu,
\end{equation}
or equivalently
\begin{equation*}\int_XR_Wf\,\mathrm{d}\nu=\int_Xf\,\mathrm{d}\nu,\mbox{ for all }f\in L^\infty(X).\end{equation*}
\end{enumerate}
\end{enumerate}

\begin{lemma}\label{lemfm1}
\begin{enumerate}
\item For measures $\mu$ on $X$ we have the implication
\textup{(\ref{eqfm5})}${}\Rightarrow{}$\textup{(\ref{eqfm4})}, but not conversely. \item
If $W$ is given and if $\nu$ and $h$ satisfy \textup{(\ref{eqfm8})} and
\textup{(\ref{eqfm7})} respectively, then
\begin{equation}\label{eqfm9}
d\mu:=h\,\mathrm{d}\nu
\end{equation}
satisfies \textup{(\ref{eqfm4})}. \item If $\mu$ satisfies \textup{(\ref{eqfm5})},
and $m_0\in L^\infty(X)$, then $S_{m_0}$ is an isometry in
$L^2(X,h\,\mathrm{d}\mu)$ if and only if
\begin{equation*}R_{m_0}h=h.\end{equation*}
\end{enumerate}
\end{lemma}

\begin{proof}
\par
(i) Suppose $\mu$ satisfies (\ref{eqfm5}). Let $f\in L^\infty(X)$. Then
\begin{equation*}\int_Xf\circ
r\,\mathrm{d}\mu=\int_X\frac{1}{\#r^{-1}(x)}\sum_{r(y)=x}f(r(y))\,\mathrm{d}\mu(x)=\int_Xf\,\mathrm{d}\mu.\end{equation*}
\par
(ii) Let $W,\nu$ and $h$ be as in the statement of part (ii) of the lemma. Then
\begin{align*}
\int_Xf\circ r\,\mathrm{d}\mu&=\int_X f\circ r\,h\,\mathrm{d}\nu=\int_XR_W(f\circ
r\,
h)\,\mathrm{d}\nu\\
&=\int_XfR_Wh\,\mathrm{d}\nu=\int_Xfh\,\mathrm{d}\nu=\int_Xf\,\mathrm{d}\mu,
\end{align*}
which is the desired conclusion (\ref{eqfm4}). It follows in particular that (\ref{eqfm5}) is strictly stronger than (\ref{eqfm4}).
\par
(iii) For $f\in L^\infty(X)$, we have
\begin{align*}
\|S_{m_0}f\|_{L^2(X,h\,\mathrm{d}\mu)}^2&=\int_X|m_0(x)f(rx)|^2h(x)\,\mathrm{d}\mu\\
&=\int_X|f(x)|^2\frac{1}{\#r^{-1}(x)}\sum_{r(y)=x}|m_0(y)|^2h(y)\,\mathrm{d}\mu(x)\\
&=\int_X|f(x)|^2R_{m_0}h(x)\,\mathrm{d}\mu(x)=\int_X|f|^2h\,\mathrm{d}\mu=\|f\|_{L^2(X,h\,\mathrm{d}\mu)}^2
\end{align*}
iff $R_{m_0}h=h$ and (iii) follows.
\end{proof}
\par
We will use standard facts from measure theory: for example, we
may identify positive Borel measures on $X$ with positive linear
functionals on $C(X)$ via
\begin{equation*}\Lambda_\omega(f)=\int_Xf\,\mathrm{d}\omega.\end{equation*}
In fact, we will identify $\Lambda_\omega$ and $\omega$. For two
measures $\mu$ and $\nu$ on $X$, we will use the notation
$\mu\prec\nu$ to denote absolute continuity. For example
$\mu\prec\nu$ holds in (\ref{eqfm9}).
\subsection{\label{Exa}Examples} We illustrate the definitions:
\begin{example}\label{exfm1}
Let $X=[0,1]=\br/\bz$. Fix $N\in\bz_+$, $N>1$. Let
\begin{equation*}r(x)=Nx\mod 1\end{equation*}
\par
Invariance:
\begin{equation}\label{eqfme4}
\int_0^1f(Nx)\,\mathrm{d}\mu(x)=\int_0^1f(x)\,\mathrm{d}\mu(x),\quad(f\in
L^\infty(\br/\bz)).
\end{equation}
\par
Strong invariance:
\begin{equation}\label{eqfme5}
\frac{1}{N}\int_0^1\sum_{k=0}^{N-1}f\left(\frac{x+k}{N}\right)\,\mathrm{d}\mu(x)=\int_0^1f(x)\,\mathrm{d}\mu(x).
\end{equation}
\par
The Lebesgue measure $\mu=\lambda$ is the unique probability
measure on $[0,1]=\br/\bz$ which satisfies (\ref{eqfme5}).
\par
Examples of measures $\mu$ on $\br/\bz$ which satisfy
(\ref{eqfme4}) but not (\ref{eqfme5}) are
\begin{itemize}
\item $\mu=\delta_0$, the Dirac mass at $x=0$; \item
$\mu=\mu_{\bf{C}}$, the Cantor middle-third measure on $[0,1]$
(see \cite{DutJo}), i.e., $\mu_{\bf{C}}$ is determined by
\begin{itemize}
\item
$\frac{1}{2}\int\left(f\left(\frac{x}{3}\right)+f\left(\frac{x+2}{3}\right)\right)\,\mathrm{d}\mu_{\bf{C}}(x)=\int
f(x)\,\mathrm{d}\mu_{\bf{C}}(x),$ \item
$\mu_{\bf{C}}([0,1])=1,$ \item
$\mu_{\bf{C}}$ is supported on the middle-third Cantor set.
\end{itemize}
\end{itemize}
\end{example}

\begin{example}\label{exfme2}
Let $X=[0,1)=\br/\bz$, $\lambda$ the Lebesgue measure, $X_{\bf C}$
the middle-third Cantor set, $\mu_{\bf C}$ the Cantor measure.
\par
$r:X\rightarrow X$, $r(x)=3x\mod 1$, $r_{\bf C}=r_{X_{\bf
C}}:X_{\bf C}\rightarrow X_{\bf C}$.
\par
Consider the following properties for a Borel probability measure
$\mu$ on $\br$:
\begin{equation}\label{eqfmeinv}
\int
f\,\mathrm{d}\mu=\frac{1}{3}\int\left(f(\frac{x}{3})+f(\frac{x+1}{3})+f(\frac{x+2}{3})\right)\,\mathrm{d}\mu(x);
\end{equation}
\begin{equation}\label{eqfmeinvc}
\int
f\,\mathrm{d}\mu=\frac{1}{2}\int\left(f(\frac{x}{3})+f(\frac{x+2}{3})\right)\,\mathrm{d}\mu(x);
\end{equation}
Then (\ref{eqfmeinv}) has a unique solution $\mu=\lambda$.
Moreover (\ref{eqfmeinvc}) has a unique solution, $\mu=\mu_{\bf
C}$, and $\mu_{\bf C}$ is supported on the Cantor set $X_{\bf C}$.
\par
Let $\br/\bz=[0,1)$. Then $\#r^{-1}(x)=3$ for all $x\in[0,1)$. If
$x=\frac{x_1}{3}+\frac{x_2}{3^2}+...$, $x_i\in\{0,1,2\}$, is the
representation of $x$ in base $3$, then $r(x)\sim(x_2,x_3,...)$,
and
$r^{-1}(x)=\{(0,x_1,x_2,...),(1,x_1,x_2,...),(2,x_1,x_2,...)\}$
\par
On the Cantor set $\#r_{\bf C}^{-1}(x)=2$ for all $x\in X_{\bf
C}$. If $x=\frac{x_1}{3}+\frac{x_2}{3^2}+...$, $x_i\in\{0,2\}$ is
the usual representation of $X_{\bf C}$ in base $3$, then
\begin{equation*}r_{\bf C}(x)=(x_2,x_3,...)\end{equation*}
and
\begin{equation*}X_{\bf C}\simeq\prod_{\mathbb{N}}\{0,2\}.\end{equation*}
\par
In the representation $\prod_{\mathbb{N}}\bz_3$ of $X=[0,1)$,
$\mu=\lambda$ is the product (Bernoulli) measure with weights
$(\frac13,\frac13,\frac13)$.
\par
In the representation $\prod_{\mathbb{N}}\{0,2\}$ of $X_{\bf C}$,
$\mu_{\bf C}$ is the product (Bernoulli) measure with weights
$(\frac12,\frac12)$.
\end{example}

\begin{example}\label{exfme3}
Let $N\in\bz_+$, $N\geq2$ and let $A=(a_{ij})_{i,j=1}^N$ be an $N$
by $N$ matrix with all $a_{ij}\in\{0,1\}$. Set
\begin{equation*}X(A):=\{(x_i)\in\prod_{\mathbb{N}}\{1,...,N\}\,|\,A(x_i,x_{i+1})=1\}\end{equation*}
and let $r=r_A$ be the restriction of the shift to $X(A)$, i.e.,
\begin{equation*}r_A(x_1,x_2,...)=(x_2,x_3,...),\quad(x=(x_1,x_2,...)\in X(A)).\end{equation*}
\begin{lemma}\label{lemfme3_1}
Let $A$ be as above. Then
\begin{equation*}\#r_A^{-1}(x)=\#\{y\in\{1,...,N\}\,|\,A(y,x_1)=1\}.\end{equation*}
\end{lemma}
\par
It follows that $r_A:X(A)\rightarrow X(A)$ is onto iff $A$ is {\it
irreducible}, i.e., iff for all $j\in\bz_N$, there exists an
$i\in\bz_N$ such that $A(i,j)=1$. \par Suppose in addition that
$A$ is {\it aperiodic}, i.e., there exists $p\in\bz_+$ such that
$A^p>0$ on $\bz_N\times\bz_N$. We have the following lemma:
\begin{lemma}[D. Ruelle, \cite{Rue89}, \cite{Bal00}]\label{lemfme3_2}
Let $A$ be irreducible and aperiodic
and let $\phi\in C(X(A))$ be given. Assume that $\phi$ is a
Lipschitz function.
\begin{enumerate}
\item
Set
\begin{equation*}(R_\phi f)(x)=\sum_{r_A(y)=x}\mathrm{e}^{\phi(y)}f(y),\mbox{ for }f\in
C(X(A)).\end{equation*} Then there exist $\lambda_0>0$,,
\begin{equation*}\lambda_0=\sup\{|\lambda|\,|\,\lambda\in\mbox{spec}(R_\phi)\},\end{equation*}
$h\in C(X(A))$ strictly positive and $\nu$ a Borel measure on
$X(A)$ such that
\begin{align*}R_\phi h&=\lambda_0 h,\\\nu R_\phi&=\lambda_0\nu,\end{align*}
and $\nu(h)=1$. The data is unique. \item In particular, setting

\begin{equation*}(R_0f)(x)=\frac{1}{\#r_A^{-1}(x)}\sum_{r_A(y)=x}f(y),\end{equation*}
we may take $\lambda_0=1$, $h=1$ and $\nu=:\mu_A$, where $\mu_A$
is a probability measure on $X(A)$ satisfying the strong
invariance property
\begin{equation*}\int_{X(A)}f\,\mathrm{d}\mu_A=\int_{X(A)}\frac{1}{\#r_A^{-1}(x)}\sum_{r_A(y)=x}f(y)\,\mathrm{d}\mu_A(x),\quad(f\in
L^\infty(X(A)).\end{equation*}
\end{enumerate}
\end{lemma}
\end{example}

\section{\label{Pos}\uppercase{Positive definite functions and dilations}} We now recall a result
relating operator systems to positive definite functions. The idea
dates back to Kolmogorov, but has been used recently in for
example \cite{FO00} and \cite{Dut3} (see also \cite{Aro}).
\begin{definition}\label{def2_1}
A map $K:X\times X\rightarrow\bc$ is called {\it positive
definite} if, for any $x_1,..,x_n\in X$ and any
$\xi_1,...,\xi_n\in\bc$,
\begin{equation*}\sum_{i,j=1}^nK(x_i,x_j)\cj\xi_i\xi_j\geq 0.\end{equation*}
\end{definition}

\begin{theorem}[Kolmogorov-Aronszajn]\label{th2_1_0}
Let $K:X\times X\rightarrow\bc$ be positive definite. Then there
exist a Hilbert space and a map $v:X\rightarrow H$ such that
\begin{gather*}\cj{\mbox{span}}\{v(x)\,|\,x\in X\}=H,\\
\ip{v(x)}{v(y)}=K(x,y),\quad(x,y\in X).\end{gather*}
Moreover $H$ and $v$ are unique up to isomorphism.
\end{theorem}
\begin{proof}
We sketch the idea of the proof. Take $H$ to be the completion of
the space
\begin{equation*}\{f:X\rightarrow\bc\,|\,f\mbox{ has finite support }\}\end{equation*}
with respect to the scalar product
\begin{equation*}\ip{f}{g}=\sum_{x,y\in X}\cj{f(x)}K(x,y)g(y).\end{equation*}
Then define $v(x):=\delta_x$.
\end{proof}

\begin{remark}\label{rem3_2}
Theorem \ref{th2_1_0} has a long history in operator theory.
The version above is purely geometric, but as noted, for example in \cite{PaSc72}
and \cite{BCR84}, it is possible to take the Hilbert space $H$ in the theorem of
the form $L^2(\Omega,\mathfrak{B}, \mu)$ where $(\Omega,\mathfrak{B}, \mu)$ is a probability space;
i.e., $\mathfrak{B}$ is a sigma-algebra on some measure space $\Omega$, $\mu$ a measure defined
on $\mathfrak{B}$, $\mu(\Omega) = 1$. In that case, $v(x, .)$ is a stochastic process. As is
well known, it is even possible to make this choice 
such that the process is
Gaussian. Examples of this include Brownian motion, and fractional Brownian
motion, see also \cite{Moh03}, \cite{Aya04}, \cite{JMR01}; -- and \cite{MoPa92} for a more
operator theoretic approach. 
\par
        For the purpose of the present discussion, it will be enough to know
the Hilbert space $H$ abstractly, but in the main part of our paper (Sections
\ref{ProLim}--\ref{Ite}), the particular function representation will be of significance. To
see this, take for example the case of the more familiar wavelet
construction from Example 1.1 above. In the present framework, the space $X$
is then the $d$-torus $\bt^d$, while the ambient dilation Hilbert space $H$ is
$L^2(\br^d)$. Since wavelet bases must be realized in the ambient Hilbert space,
it is significant to have much more detail than is encoded in the purely
geometric data of Theorem \ref{th2_1_0}. Even when comparing with the function
theoretic version of \cite{PaSc72}, the wavelet example illustrates that it is
significant to go beyond probability spaces.
\par
       One of our aims is to offer a framework for more general wavelet
bases, including state spaces in symbolic dynamics and Julia sets (such as
\cite{BrTo05}.) A main reason for the usefulness of wavelet bases is their
computational features. As is well known \cite{San59}, there are many function
theoretic orthonormal bases (ONB), or Parseval frames in analysis where the
basis coefficients do not lend themselves practical algorithmic schemes. If
for example we are in $L^2(\br^d)$, then the computation of each basis
coefficients typically involves a separate integration over $\br^d$; not at all
a computationally attractive proposition. 
\par
       What our present approach does is that it selects a subspace of the
ambient Hilbert space which is computationally much more feasible. As
stressed in \cite{BJMP05} and \cite{Jor06}, such a selection corresponds to a choice
of resolution, a notion from optics; and one dictated in turn by
applications. In the present setup, the chosen resolution corresponds to an
initial space, which in this context may be encoded by $X$ from Theorem \ref{th2_1_0}
above. As we will see later, there are ways to do this in such that the
computation of basis coefficients becomes algorithmic. We will talk about
wavelet bases in this much more general contest, even though wavelets are
traditionally considered only in $L^2(\br^d)$. With good choices, we find that
computation of the corresponding basis coefficients may be carried with a
certain recursive algorithms involving only matrix iteration; much like in
the familiar case of Gram-Schmidt algorithms.
\end{remark}
\begin{theorem}\label{th2_2}
Let $K$ be a positive definite map on a set $X$. Let
$s:X\rightarrow X$ be a map that is compatible with $K$ in the
sense that
\begin{equation}\label{eq2_1}
K(s(x),s(y))=K(x,y),\quad(x,y\in X).
\end{equation}
Then there exists a Hilbert space $H$, a map $v:X\rightarrow H$
and a unitary operator $U$ on $H$ such that
\begin{equation}\label{eq2_2}
\ip{v(x)}{v(y)}=K(x,y),\quad(x,y\in X),
\end{equation}
\begin{equation}\label{eq2_3}
\cj{\mbox{span}}\{U^{-n}(v(x))\,|\,x\in X,n\geq 0\}=H,
\end{equation}
\begin{equation}\label{eq2_4}
Uv(x)=v(s(x)),\quad(x\in X).
\end{equation}
Moreover, this is unique up to an intertwining isomorphism.
\end{theorem}

\begin{proof}
Let $\tilde X:=X\times\bz$. Define $\tilde K:\tilde X\times\tilde
X\rightarrow\bc$ by
\begin{equation*}\tilde K((x,n),(y,m))=K(s^{n+M}(x),s^{m+M}(y)),\quad(x,y\in
X,n,m\in\bz),\end{equation*} where $M\geq\max\{-m,-n\}$.
\par
The compatibility condition (\ref{eq2_1}) implies that the
definition does not depend on the choice of $M$. We check that
$\tilde K$ is positive definite. Take $(x_i,n_i)\in\tilde X$ and
$\xi_i\in\bc$. Then, for $M$ big enough we have:
\begin{equation*}\sum_{i,j}\tilde
K((x_i,n_i),(x_j,n_j))\cj\xi_i\xi_j=\sum_{i,j}K(s^{M+n_i}(x_i),s^{M+n_j}(x_j))\cj\xi_i\xi_j\geq0.\end{equation*}
Using now the Kolmogorov construction (see Theorem \ref{th2_1_0}),
there exists a Hilbert space $H$ and a map $\tilde v:\tilde
X\rightarrow H$ such that
\begin{gather*}\ip{\tilde v(x,n)}{\tilde v(y,m)}=\tilde
K((x,n),(y,m)),\quad((x,m),(y,n)\in\tilde X),\\
\cj{\mbox{span}}\{\tilde v(x,m)\,|\,(x,m)\in\tilde X\}=H.\end{gather*}
Define $v$ by
\begin{equation*}v(x)=\tilde v(x,0),\quad(x\in X).\end{equation*}
Then (\ref{eq2_2}) is satisfied. Define
\begin{equation*}U\tilde v(x,n)=\tilde v(x,n+1),\quad((x_n)\in\tilde X).\end{equation*}
$U$ is well defined and an isometry because, for $M$ sufficiently
big,
\begin{align*}\ip{\tilde v(x,n+1)}{\tilde v(y,m+1)}&=K(s^{M+n+1}(x),s^{M+m+1}(y))\\
&=K(s^{M+n}(x),s^{M+m}(y))=\ip{\tilde v(x,n)}{\tilde v(y,m)}.\end{align*}
$U$ has dense range  so $U$ is unitary. Also (\ref{eq2_3}) is
immediate (we need only $n\geq 0$ because $U^n(v(x))=v(s^n(x))$,
for $n\geq0$, will follow form (\ref{eq2_4})).
\par
For (\ref{eq2_4}) we compute
\begin{align*}\ip{Uv(x)}{\tilde v(y,n)}&=\tilde
K((x,1),(y,n))=K(s^{M+1}(x),s^{M+n}(y))\\
&=K(s^M(s(x)),s^{M+n}(y))=\ip{v(s(x))}{\tilde v(y,n)}.\end{align*}
For uniqueness, if $H',v', U'$ satisfy the same conditions, then the formula\\
$W(U^nv(x))=U'^nv'(x)$ defines an intertwining isomorphism.
\end{proof}

\begin{theorem}\label{th2_3}
Let $\mathcal{A}$ be a unital $C^*$-algebra, $\alpha$ an
endomorphism on $\mathcal{A}$, $\mu$ a state on $\mathcal{A}$ and,
$m_0\in\mathcal{A}$, such that
\begin{equation}\label{eq2_5}
\mu(m_0^*\alpha(f)m_0)=\mu(f),\quad(f\in\mathcal{A}).
\end{equation}
Then there exists a Hilbert space $H$, a representation $\pi$ of
$\mathcal{A}$ on $H$, $U$ a unitary on $H$, and a vector
$\varphi\in\mathcal{A}$, with the following properties:
\begin{equation}\label{eq2_6}
U\pi(f)U^*=\pi(\alpha(f)),\quad(f\in\mathcal{A}),
\end{equation}
\begin{equation}\label{eq2_7}
\ip{\varphi}{\pi(f)\varphi}=\mu(f),\quad(f\in\mathcal{A}),
\end{equation}
\begin{equation}\label{eq2_8}
U\varphi=\pi(\alpha(1)m_0)\varphi
\end{equation}
\begin{equation}\label{eq2_9}
\cj{\mbox{span}}\{U^{-n}\pi(f)\varphi\,|\,n\geq0,f\in\mathcal{A}\}=H.
\end{equation}
Moreover, this is unique up to an intertwining isomorphism.
\par
We call $(H,U,\pi,\varphi)$ the covariant system associated to
$\mu$ and $m_0$.
\end{theorem}

\begin{proof}
Define $K$ and $s$ by
\begin{equation*}K(x,y)=\mu(x^*y), \quad
s(x)=\alpha(x)m_0,\quad(x,y\in\mathcal{A}).\end{equation*} $K$ is positive
definite and compatible with $s$ so, with Theorem \ref{th2_2},
there exists a Hilbert space $H$, a map $v$ from $\mathcal{A}$ to
$H$, and a unitary $U$ with the mentioned properties.
\par
Define $\varphi=v(1)$,
\begin{equation*}\pi(f)(U^{-n}v(x))=U^{-n}v(\alpha^n(f)x),\quad(f,x\in\mathcal{A},n\geq0).\end{equation*}
Some straightforward computations show that $\pi$ is a well
defined representation of $\mathcal{A}$ that satisfies all
requirements.
\end{proof}

\begin{corollary}\label{cor2_4}
Let $X$ be a measure space, $r:X\rightarrow X$ a measurable, onto
map and $\mu$ a probability measure on $X$ such that
\begin{equation}
\int_Xf\,\mathrm{d}\mu=\int_X\frac{1}{\#r^{-1}(x)}\sum_{r(y)=x}f(y)\,\mathrm{d}\mu(x).
\end{equation}
Let $h\in L^1(X)$, $h\geq0$ such that \begin{equation*}\frac{1}{\#
r^{-1}(x)}\sum_{r(y)=x}|m_0(y)|^2h(y)=h(x),\quad(x\in X).\end{equation*} Then
there exists (uniquely up to isomorphisms) a Hilbert space $H$, a
unitary $U$, a representation $\pi$ of $L^\infty(X)$ and a vector
$\varphi\in H$ such that
\begin{gather*}U\pi(f)U^{-1}=\pi(f\circ r),\quad(f\in L^\infty(X)),\\
\ip{\varphi}{\pi(f)\varphi}=\int_Xfh\,\mathrm{d}\mu,\quad(f\in
L^\infty(X)),\\
U\varphi=\pi(m_0)\varphi,\vphantom{\int_X}\\
\cj{\mbox{span}}\{U^{-n}\pi(f)\varphi\,|\,n\geq0,f\in
L^\infty(X)\}=H.\end{gather*} We call $(H,U,\pi,\varphi)$ the covariant
system associated to $m_0$ and $h$.
\end{corollary}

\begin{proof}
Take $\mu(f)=\int_Xfh\,\mathrm{d}\mu$, $\alpha(f)=f\circ r$; and use
Theorem \ref{th2_3}.
\end{proof}
\par
We regard Theorem \ref{th2_3} as a dilation result. In this
context we have a second closely related result:
\begin{theorem}\label{th2_5}
\textup{(i)} Let $H$ be a Hilbert space, $S$ an isometry on $H$. Then there
exist a Hilbert space $\hat H$ containing $H$ and a unitary $\hat
S$ on $\hat H$ such that
\begin{equation}\label{eq2_5_1}
\hat S|_{H}=S,
\end{equation}
\begin{equation}\label{eq2_5_2}
\cj{\bigcup_{n\geq0}\hat S^{-n}H}=\hat H.
\end{equation}
Moreover these are unique up to an intertwining isomorphism. \par
\textup{(ii)} If $\mathcal{A}$ is a $C^*$-algebra, $\alpha$ is an
endomorphism on $\mathcal{A}$ and $\pi$ is a representation of
$\mathcal{A}$ on $H$ such that
\begin{equation}\label{eq2_5_3}
S\pi(g)=\pi(\alpha(g))S,\quad(g\in\mathcal{A});
\end{equation}
then there exists a unique representation $\hat\pi$ on $\hat H$
such that
\begin{equation}\label{eq2_5_4}
\hat\pi(g)|_H=\pi(g),\quad(g\in\mathcal{A}),
\end{equation}
\begin{equation}\label{eq2_5_5}
\hat S\hat\pi(g)=\hat\pi(\alpha(g))\hat S,\quad(g\in\mathcal{A}).
\end{equation}
\end{theorem}

\begin{proof}
(i) Consider the set of symbols
\begin{equation*}\mathcal{H}_{sym}:=\{\sum_{j\in\bz}S^j\xi_j\,|\,\xi_j\in
H,\xi_j=0\mbox{ except for finiteley many }j\mbox{'s}\}.\end{equation*} Define
the scalar product
\begin{equation}\label{diez}
\Bip{\sum_{i\in\bz}S^i\xi_i}{\sum_{j\in\bz}S^j\eta_j}=\sum_{i,j\in\bz}\bip{S^{i+m}\xi_i}{S^{j+m}\eta_j},
\end{equation}
where $m$ is chosen sufficiently large, such that $i+m,j+m\geq0$
for all $i,j\in\bz$ with $\xi_i\neq0,\eta_j\neq0$.
\par
Since $S$ is an isometry this definition does not depend on the
choice of $m$. We denote the completion of $\mathcal{H}_{sym}$
with this scalar product by $\hat H$. $H$ can be isometrically
identified with a subspace of $\hat H$ by
\begin{equation*}\xi\mapsto\sum_{i\in\bz}S^i\xi_i,\mbox{ where
}\xi_i=\left\{\begin{array}{ccc}0&\mbox{if}&i\neq0\\
\xi&\mbox{if}&i=0.\end{array}\right.\end{equation*} Define
\begin{equation*}\hat S(\sum_{i\in\bz}S^i\xi_i)=\sum_{i\in\bz}S^{i+1}\xi_i.\end{equation*}
\par
In the definition of $\hat H$, we use (\ref{diez}) as an inner
product, and we set
\begin{equation*}\hat H=\left(\mathcal{H}_{sym}/\{\sum_jS^j\xi_j\,|\,\sum_{i,j}\bip{S^{i+m}\xi_i}{S^{j+m}\xi_j}=0\}\right)^{\wedge}\end{equation*}
where $\wedge$ stands for completion.
\par
Since $\xi=S^{-1}(S\xi)$ in $\mathcal{H}_{sym}$, for $\xi\in H$,
we get natural isometric embeddings as follows, see
(\ref{eq2_5_2}),
\begin{equation*}H\subset\hat S^{-1}H\subset\hat S^{-2}H\subset...\subset\hat
S^{-n}H\subset\hat S^{-n-1}H\subset...\end{equation*}
\par
It can be checked that $\hat H$ and $\hat S$ satisfy the
requirements.
\par
(ii) We know that the spaces
\begin{equation*}\{\hat S^{-n}\xi\,|\,n\geq 0,\xi\in H\}\end{equation*}
span a dense subspace of $\hat H$. Define
\begin{equation*}\hat\pi(g)(\hat S^{-n}\xi)=\hat
S^{-n}\pi(\alpha^n(g))\xi,\quad(g\in\mathcal{A},n\geq0,\xi\in
H).\end{equation*} We check only that $\hat\pi(g)$ is a well defined, bounded
operator, the rest of our claims follow from some elementary
computations. Take $m$ large:
\begin{align*}
\|\hat\pi(g)(\sum_i\hat S^{-n_i}\xi_i)\|^2&=\|\sum_i\hat
S^{-n_i}\pi(\alpha^{n_i}(g))\xi_i\|^2\\
&=\|\sum_{i}\hat S^{m-n_i}\pi(\alpha^{n_i}(g))\xi_i\|^2\\
&=\|\sum_i\pi(\alpha^m(g))\hat S^{m-n_i}\xi_i\|^2\\
&\leq\|g\|^2\|\sum_i\hat S^{-n_i}\xi_i\|^2.\qedhere
\end{align*}
\end{proof}

\begin{example}\label{exbaggett}
This example is from \cite{BCMO}, and it illustrates the
conclusions in Theorem \ref{th2_5}.
\par

Consider\\
1. $H=l^2(\mathbb{N}_0),$\\
2. $S(c_0,c_1,..)=(c_1,c_2,...),$ the unilateral shift,\\
3. $\delta_k(j)=\delta_{k,j}=$ Kronecker delta, for
$k,j\in\mathbb{N}_0$,\\
4. $\pi(g_k)\delta_j:=\exp(i2\pi k2^{-j})\delta_j,\quad
j\in\mathbb{N}_0,k\in\bz[1/2]$.

\par
When Theorem \ref{th2_5} is applied we get:\\
1'. The dilation Hilbert space $\hat H$ is $l^2(\bz)$. \\
2'. $\hat S$ is the bilateral shift on $l^2(\bz)$ i.e., $\hat
S\delta_j=\delta_{j-1}$ for $j\in\bz$.\\
3'. Same as in 3. but for $k,j\in\bz$.\\
4'. The operator $\hat\pi(g_k)$ is given by the same formula 4.,
but for $j\in\bz$.
\par
The commutation relation (\ref{eq2_5_5}) now takes the form
\begin{equation}\label{eqbaggett5}
\hat S\hat\pi(g_k)=\hat\pi(g_{2^{-1}k})\hat S\mbox{ on
}l^2(\bz),\mbox{ for }k\in\bz[1/2];
\end{equation}
and
\begin{equation}\label{eqbaggett6}
\hat S^{-n}
H=\cj{\mbox{span}}\{\delta_{-n},\delta_{-n+1},\delta_{-n+2},...\}\subset\hat
H.
\end{equation}

\end{example}

\subsection{\label{Ope}Operator valued filters}
\par
In this subsection we study the multiplicity configurations of the
representations $\pi$ from above. Our first result shows that the
two functions $m_0$, and $h$ in Section \ref{Tra} may be operator
valued. The explicit multiplicity functions are then calculated in
the next section.
\begin{corollary}\label{cor2_6}
Let $X,r,$ and $\mu$ be as in Corollary \ref{cor2_4}. Let $I$ be a
finite or countable set. Suppose
$H:X\rightarrow\mathcal{B}(l^2(I))$ has the property that
$H(x)\geq 0$ for almost every $x\in X$, and $H_{ij}\in L^1(X)$ for
all $i,j\in I$. Let $M_0:X\rightarrow\mathcal{B}(l^2(I))$ such
that $x\mapsto\|M_0(x)\|$ is essentially bounded. Assume in
addition that
\begin{equation}\label{eq2_6_1}
\frac{1}{\#r^{-1}(x)}\sum_{r(y)=x}M_0^*(y)H(y)M_0(y)=H(x),\mbox{ for a.e. }x\in X.
\end{equation}
Then there exists a Hilbert space $\hat K$, a unitary operator
$\hat U$ on $\hat K$, a representation $\hat\pi$ of $L^\infty(X)$
on $\hat K$, and a family of vectors $(\varphi_i)\in\hat K$, such
that:
\begin{gather*}\hat U\hat\pi(g)\hat U^{-1}=\hat\pi(g\circ r),\quad(g\in L^\infty(X)),\\
\hat U\varphi_i=\sum_{j\in I}\hat\pi((M_0)_{ji})\varphi_j,\quad(i\in I),\vphantom{\int_X}\\
\ip{\varphi_i}{\hat\pi(f)\varphi_j}=\int_XfH_{ij}\,\mathrm{d}\mu,\quad(i,j\in I,f\in L^\infty(X)),\\
\cj{\mbox{span}}\{\hat\pi(f)\varphi_i\,|\,n\geq0,f\in L^\infty(X),i\in I\}=\hat K.\end{gather*}
These are unique up to an intertwining unitary isomorphism.
(All functions are assumed weakly measurable in the sense that $x\mapsto\ip{\xi}{F(x)\eta}$ is measurable for all $\xi,\eta\in l^2(I)$.)
\end{corollary}
\begin{proof}
Consider the Hilbert space
\begin{equation*}K:=\{f:X\rightarrow\bc^I\,|\,f\mbox{ is measurable }, \int_X\ip{f(x)}{H(x)f(x)}\,\mathrm{d}\mu(x)<\infty\}.\end{equation*}
Define $S$ on $K$ by
\begin{equation*}(Sf)(x)=M_0(x)(f(r(x))),\quad(x\in X,f\in K).\end{equation*}
We check that $S$ is an isometry. For $f,g\in K$:
\begin{align*}
\ip{Sg}{Sf}&=\int_X\ip{M_0(x)g(r(x))}{H(x)M_0(x)f(r(x))}\,\mathrm{d}\mu(x)\\
&=\int_X\ip{g(r(x))}{M_0(x)^*H(x)M_0(x)f(r(x))}\,\mathrm{d}\mu(x)\\
&=\int_X\frac{1}{\#r^{-1}(x)}\sum_{r(y)=x}\ip{g(x)}{M_0(y)^*H(y)M_0(y)f(x)}\,\mathrm{d}\mu(x)\\
&=\int_X\ip{g(x)}{H(x)f(x)}\,\mathrm{d}\mu(x)=\ip{g}{f},
\end{align*}
where we used (\ref{eq2_6_1}) in the last step. The converse
implication holds as well, i.e., if $S$ is an isometry then
(\ref{eq2_6_1}) is satisfied.\par
 Define now
\begin{equation*}(\pi(g)f)(x)=g(x)f(x),\quad(x\in X,g\in L^\infty(X),f\in K).\end{equation*}
$\pi$ defines a representation of $L^\infty(X)$ on $K$. Moreover, the covariance relation is satisfied
\begin{equation*}S\pi(g)=\pi(g\circ r)S.\end{equation*}
Then we use Theorem \ref{th2_5} to obtain a Hilbert space $\hat K$
containing $K$, a unitary $\hat U:=\hat S$ on $\hat K$, and a
representation $\hat\pi$ on $\hat K$ that dilate $S$ and $\pi$.
\par
Define $\varphi_i\in K\subset\hat K$,
\begin{equation*}\varphi_i(x):=\delta_i,\mbox{ for all }x\in X,\quad(i\in I).\end{equation*}
We have that
\begin{gather*}\ip{\varphi_i}{\hat\pi(f)\varphi_j}=\int_X\ip{\delta_i}{H(x)(f(x)\delta_j)}\,\mathrm{d}\mu(x)=\int_Xf(x)H_{ij}(x)\,\mathrm{d}\mu(x),\\
(\hat U\varphi_i)(x)=(S\varphi_i)(x)=M_0(x)\delta_i=((M_0)_{ji}(x))_{j\in I}=(\sum_{j\in I}\hat\pi((M_0)_{ji})\varphi_j)(x).\end{gather*}
Also it is clear that
\begin{equation*}\cj{\mbox{span}}\{\hat\pi(f)\varphi_i\,|\,f\in L^\infty(X),i\in I\}=K.\end{equation*}
These relations, together with Theorem \ref{th2_5}, prove our
assertions.
\end{proof}

\section{\label{Mul}\uppercase{Multiplicity theory}}
\par
 One of the tools from operator theory which has been especially useful in the analysis of wavelets is multiplicity theory
for abelian $C^*$-algebras $\mathcal{A}$.

We first recall a few well known facts, see e.g., \cite{N}. By
Gelfand's theorem, every abelian $C^*$-algebra with unit is $C(X)$
for a compact Hausdorff space X; and every representation of
$\mathcal{A}$ is the orthogonal sum of cyclic representations.
While the cardinality of the set of cyclic components in this
decomposition is an invariant, the explicit determination of the
cyclic components is problematic, as the construction depends on
Zorn's lemma.
   So for this reason, it is desirable to turn the abstract spectral theorem for representations into a concrete one. In the
concrete spectral representation, $C(X)$ is represented as an
algebra of multiplication operators on a suitable $L^2$-space; as
opposed to merely an abstract Hilbert space. When we further
restrict attention to normal representations of $\mathcal{A}$, we
will be working with the algebra $L^\infty(X)$ defined relative to
the Borel sigma-algebra of subsets in $X$.
\par
With this, we are able to compute a concrete spectral
representation, and thereby to strengthen the conclusion from
Theorem \ref{th2_5}.
\par
 Our $L^2$-space which carries the representation may be realized concretely when the additional structure from Section \ref{Tra}
is introduced, i.e., is added to the assumptions in Theorem
\ref{th2_5}. Hence, we will work with the given finite-to-one
mapping $r:X\rightarrow X$, and the measure $\mu$ from before.
Recall from Section \ref{Fun} that $\mu$ is assumed strongly
$r$-invariant.
\par
Theorem \ref{th2_5} provides an abstract unitary dilation of a
given covariant system involving a representation $\pi$ and a
fixed isometry $S$ on a Hilbert space $H$. In the present section,
we specialize the representation $\pi$ in Theorem \ref{th2_5} to
the algebra $\mathcal{A}=L^\infty(X)$, and
$\alpha:\mathcal{A}\rightarrow\mathcal{A}$, is $\alpha(g):=g\circ
r$. \par While our conclusion from Theorem \ref{th2_5} still
offers a unitary dilation $U$ in an abstract Hilbert space $\hat
H$, we are now able to show that $\hat H$ has a concrete spectral
representation. Since $\hat H$ is the closure of an ascending
union of resolution subspaces defined from $U$, the question
arises as to how the multiplicities of the restricted
representations of the resolution subspaces in $\hat H$ are
related to one-another.
\par
The answer to this is known in the case of wavelets, see e.g.,
\cite{BM}. In this section we show that there is a version of the
Baggett et al multiplicity formula in the much more general
setting of Theorem \ref{th2_5}. In particular, we get the
multiplicity formula in the applications where $X$ is a Julia set,
or a state space of sub-shift dynamical system. As we noted in
Section \ref{Fun} above, each of these examples carries a natural mapping
$r$, and a strongly $r$-invariant measure $\mu$.
\par
Consider $X$ a measure space, $r:X\rightarrow X$ an onto,
measurable map such that $\#r^{-1}(x)<\infty$ for all $x\in X$.
Let $\mu$ be a measure on $X$ such that
\begin{equation}\label{eqmul1}
\int_Xf\,\mathrm{d}\mu=\int_X\frac{1}{\#r^{-1}(x)}\sum_{r(y)=x}f(y)\,\mathrm{d}\mu(x),\quad(f\in
L^\infty(X)).
\end{equation}
\par
Suppose now that $H$ is a Hilbert space with an isometry $S$ on it
and with a normal representation $\pi$ of $L^\infty(X)$ on $H$
that satisfies the covariance relation
\begin{equation}\label{eqmul2}
S\pi(g)=\pi(g\circ r)S,\quad(g\in L^\infty(X)).
\end{equation}
\par
Theorem \ref{th2_5} shows that there exists a Hilbert space $\hat
H$ containing $H$, a unitary $\hat S$ on $\hat H$  and a
representation $\hat\pi$ of $L^\infty(X)$ on $\hat H$ such that:
\begin{gather*}(V_n:=\hat S^{-n}(H))_n\mbox{ form an increasing sequence of
subspaces with dense union},\\
\hat S|_H=S,\\
\hat\pi|_H=\pi,\\
\hat S\hat\pi(g)=\hat\pi(g\circ r)\hat S.\end{gather*}
\begin{theorem}\label{thmul1}
\textup{(i)} $V_1=\hat S^{-1}(H)$ is invariant for the representation
$\hat\pi$. The multiplicity functions of the representation
$\hat\pi$ on $V_1$, and on $V_0=H$, are related by
\begin{equation}\label{eqmul3}
m_{V_1}(x)=\sum_{r(y)=x}m_{V_0}(y),\quad(x\in X).
\end{equation}
\par
\textup{(ii)} If $W_0:=V_1\ominus V_0=\hat S^{-1}H\ominus H$, then
\begin{equation}\label{eqmul4}
m_{V_0}(x)+m_{W_0}(x)=\sum_{r(y)=x}m_{V_0}(y),\quad(x\in X).
\end{equation}
\end{theorem}
\begin{proof}
Note that $\hat S$ maps $V_1$ to $V_0$, and the covariance
relation implies that the representation $\hat\pi$ on $V_1$ is
isomorphic to the representation $\pi^r:g\mapsto\pi(g\circ r)$ on
$V_0$. Therefore we have to compute the multiplicity of the
latter, which we denote by $m^r_{V_0}$.
\par
By the spectral theorem there exists a unitary isomorphism
$J:H(=V_0)\rightarrow L^2(X,m_{V_0},\mu)$, where, for a {\it
multiplicity function} $m:X\rightarrow\{0,1,...,\infty\}$, we use
the notation:
\begin{equation*}
L^2(X,m,\mu):=\{f:X\rightarrow\cup_{x\in X}\bc^{m(x)}\,|\,
f(x)\in\bc^{m(x)},\int_X\|f(x)\|^2\,\mathrm{d}\mu(x)<\infty\}.\end{equation*} In
addition $J$ intertwines $\pi$ with the representation of
$L^\infty(X)$ by multiplication operators, i.e.,
\begin{equation*}(J\pi(g)J^{-1}(f))(x)=g(x)f(x)\,\quad(g\in L^\infty(X),f\in
L^2(X,m_{V_0},\mu),x\in X).\end{equation*}

\begin{remark}
Here we are identifying $H$ with $L^2(X,m_{V_0},\mu)$ via the {\it
spectral representation}. We recall the details of this
representation $H\ni f\mapsto\tilde f\in L^2(X,m_{V_0},\mu)$.
\par
Recall that any normal representation $\pi\in Rep(L^\infty(X),H)$
is the orthogonal sum
\begin{equation}\label{eqmulrem*}
H=\sum_{k\in C}^{ }{ }^\oplus[\pi(L^\infty(X))k],
\end{equation}
where the set $C$ of vectors $k\in H$ is chosen such that
\begin{itemize}
\item $\quad\quad\quad\quad\quad\quad\quad\quad\quad\quad\|k\|=1,$
\begin{equation}\label{eqmulrem**}
\ip{k}{\pi(g)k}=\int_Xg(x)v_k(x)^2\,\mathrm{d}\mu(x),\mbox{ for all }k\in
C;
\end{equation}
\item $\ip{k'}{\pi(g)k}=0,\quad g\in L^\infty(X),k,k'\in C,k\neq
k';\mbox{ orthogonality}.$ \end{itemize}
\par
The formula (\ref{eqmulrem*}) is obtained by a use of Zorn's
lemma. Here, $v_k^2$ is the Radon-Nikodym derivative of
$\ip{k}{\pi(\cdot)k}$ with respect to $\mu$, and we use that $\pi$
is assumed normal.
\par
For $f\in H$, set
\begin{equation*}f=\sum_{k\in C}^{ }{ }^\oplus\pi(g_k)k,\quad g_k\in L^\infty(X)\end{equation*}
and
\begin{equation*}\tilde f=\sum_{k\in C}^{ }{ }^\oplus g_kv_k\in L_\mu^2(X,l^2(C)).\end{equation*}
Then $Wf=\tilde f$ is the desired spectral transform, i.e.,
\begin{gather*}W\mbox{ is unitary},\\
W\pi(g)=M(g)W,\\
\intertext{and}\|\tilde f(x)\|^2=\sum_{k\in C}|g_k(x)v_k(x)|^2.\end{gather*}
Indeed, we have
\begin{align*}\int_X\|\tilde f(x)\|^2\,\mathrm{d}\mu(x)&=\int_X\sum_{k\in
C}|g_k(x)|^2v_k(x)^2\,\mathrm{d}\mu(x)=\sum_{k\in
C}\int_X|g_k|^2v_k^2\,\mathrm{d}\mu\\
&=\sum_{k\in C}\bip{k}{\pi(|g_k|^2)k}=\sum_{k\in
C}\|\pi(g_k)k\|^2=\|\sum_{k\in C}^{ }{
}^\oplus\pi(g_k)k\|^2_H\\&=\|f\|^2_{H}.\end{align*}
\par
It follows in particular that the multiplicity function
$m(x)=m_{H}(x)$ is
\begin{equation*}m(x)=\#\{k\in C\,|\,v_k(x)\neq0\}.\end{equation*}
Setting
\begin{equation*}X_i:=\{x\in X\,|\,m(x)\geq i\},\quad (i\geq 1),\end{equation*}
we see that
\begin{equation*}H\simeq\sum^{}{}^\oplus L^2(X_i,\mu)\simeq L^2(X,m,\mu),\end{equation*}
and the isomorphism intertwines $\pi(g)$ with multiplication
operators.
\end{remark}
\par
Returning to the proof of the theorem, we have to find the similar
form for the representation $\pi^r$. Let
\begin{equation}\label{eqmul5}
\tilde m(x):=\sum_{r(y)=x}m_{V_0}(y),\quad(x\in X). \end{equation}
Define the following unitary isomorphism:
\begin{gather*}L:L^2(X,m_{V_0},\mu)\rightarrow L^2(X,\tilde m,\mu),\\
(L\xi)(x)=\frac{1}{\sqrt{\#r^{-1}(x)}}(\xi(y))_{r(y)=x}.\end{gather*}
(Note that the dimensions of the vectors match because of
(\ref{eqmul5})). This operator $L$ is unitary. For $\xi\in
L^2(X,m_{V_0},\mu)$, we have
\begin{align*}
\|L\xi\|^2_{L^2(X,m_{V_0},\mu)}&=\int_X\|L\xi(x)\|^2\,\mathrm{d}\mu(x)\\
&=\int_X\frac{1}{\#r^{-1}(x)}\sum_{r(y)=x}\|\xi(y)\|^2\,\mathrm{d}\mu(x)\\
&=\int_X\|\xi(x)\|^2\,\mathrm{d}\mu(x).
\end{align*}

And $L$ intertwines the representations. Indeed, for $g\in
L^\infty(X)$,
\begin{equation*}L(g\circ r\,\xi)(x)=(g(r(y))\xi(y))_{r(y)=x}=g(x)L(\xi)(x).\end{equation*}
Therefore, the multiplicity of the representation
$\pi^r:g\mapsto\pi(g\circ r)$ on $V_0$ is $\tilde m$, and this
proves (i).
\par
(ii) follows from (i).
\par
{\it Conclusions.} By definition, if $k\in C$,
\begin{gather*}\ip{k}{\pi(g)k}=\int_Xg(x)v_k(x)^2\,\mathrm{d}\mu(x),\mbox{\quad and }\\
\ip{k}{\pi^r(g)k}=\int_Xg(r(x))v_k(x)^2\,\mathrm{d}\mu(x)=\int_Xg(x)\frac{1}{\#r^{-1}(x)}\sum_{r(y)=x}v_k(x)^2\,\mathrm{d}\mu(x);\end{gather*}
and so
\begin{align*}
m^r(x)&=\#\{k\in
C\,|\,\sum_{r(y)=x}v_k(y)^2>0\}\\
&=\sum_{r(y)=x}\#\{k\in C\,|\, v_k(y)^2>0\}\\
&=\sum_{r(y)=x}m(y).
\end{align*}
Let $C^m(x):=\{k\in C\,|\,v_k(x)\neq0\}$. Then we showed that
\begin{equation*}C^m(x)=\bigcup_{y\in X, r(y)=x}C^m(y)\end{equation*}
and that $C^m(y)\cap C^m(y')=\emptyset$ when $y\neq y'$ and
$r(y)=r(y')=x.$ Setting $\mathcal{H}(x)=l^2(C^m(x))$, we have
\begin{equation*}\mathcal{H}(x)=l^2(C^m(x))=\sum_{r(y)=x}^{}{}^\oplus
l^2(C^m(y))=\sum_{r(y)=x}^{}{}^\oplus\mathcal{H}(y).\qedhere\end{equation*}
\end{proof}
\begin{remark}
There are many representations $(\pi,U,\hat H)$ for which
\begin{equation*}U\pi(g)U^{-1}=\pi(g\circ r),\quad(g\in C(X)),\end{equation*} holds; but for which
the spectral measures of $\pi$ are not absolutely continuous;
i.e., the measure
\begin{equation*}g\mapsto\bip{\hat h}{\pi(g)\hat h}=\int_Xg(x)\,\mathrm{d}\mu_{\hat h}(x)\end{equation*}
is singular with respect to the Julia-measure $\mu$ for some $\hat
h\in\hat H$. But for the purpose of wavelet analysis, it is
necessary to restrict our attention normal representations $\pi$.
\end{remark}

\section{\label{ProLim}\uppercase{Projective limits}}
We work in either the category of measure spaces or topological
spaces. \begin{definition} Let $r:X\rightarrow X$ be onto, and
assume that $\# r^{-1}(x)<\infty$ for all $x\in X$. We define the
{\it projective limit} of the system:
\begin{equation}\label{eqpr_1} X\stackrel{r}{\leftarrow}
X\stackrel{r}{\leftarrow}X\stackrel{r}{\leftarrow}...X_\infty
\end{equation} as
\begin{equation*}X_\infty:=\{\hat x=(x_0,x_1,...)\,|\,r(x_{n+1})=x_n,\mbox{ for all }n\geq0\}\end{equation*}

\end{definition}
\par Let $\theta_n:X_\infty\rightarrow X$ be the projection
onto the $n$-th component:
\begin{equation*}\theta_n(x_0,x_1,...)=x_n,\quad((x_0,x_1,...)\in X_\infty).\end{equation*}
Taking inverse images of sets in $X$ through these projections, we
obtain a sigma algebra on $X_\infty$, or a topology on $X_\infty$.
\par We have an induced mapping $\hat r: X_\infty\rightarrow
X_\infty$ defined by
\begin{equation}\label{eqpr_3}
\hat r(\hat x)=(r(x_0),x_0,x_1,...),\mbox{ and with inverse }\hat
r^{-1}(\hat x)=(x_1,x_2,...). \end{equation} so $\hat r$ is an
automorphism, i.e., $\hat r\circ\hat r^{-1}=\mbox{id}_{X_\infty}$
and $\hat r^{-1}\circ\hat r=\mbox{id}_{X_\infty}$.
\par
Note that
\begin{equation*}\theta_n\circ\hat r=r\circ \theta_n=\theta_{n-1}.\end{equation*}
\begin{equation*}
\setlength{\unitlength}{0.075\textwidth}
\begin{picture}(4.8,3.5)(-0.9,-0.75)
\put(0.075,2.25){\vector(1,0){2.85}}
\put(0,2){\vector(3,-2){3}}
\put(3.375,1.925){\vector(0,-1){1.85}}
\put(-0.375,2.25){\makebox(0,0){$X_\infty$}}
\put(3.375,2.25){\makebox(0,0){$X$}}
\put(3.375,-0.25){\makebox(0,0){$X$}}
\put(1.5,2.375){\makebox(0,0)[b]{$\theta_n$}}
\put(1.5,1){\makebox(0,0)[tr]{$\theta_{n-1}$}}
\put(3.5,1){\makebox(0,0)[l]{$r$}}
\end{picture}
\begin{picture}(4.8,3.5)(-0.9,-0.75)
\put(0.075,2.25){\vector(1,0){2.85}}
\put(0,2){\vector(3,-2){3}}
\put(3.375,1.925){\vector(0,-1){1.85}}
\put(-0.375,2.25){\makebox(0,0){$X_\infty$}}
\put(3.375,2.25){\makebox(0,0){$X_\infty$}}
\put(3.375,-0.25){\makebox(0,0){$X$}}
\put(1.5,2.375){\makebox(0,0)[b]{$\hat r$}}
\put(1.5,1){\makebox(0,0)[tr]{$\theta_{n-1}$}}
\put(3.5,1){\makebox(0,0)[l]{$\theta_n$}}
\end{picture}
\begin{picture}(3.7,3.5)(-0.85,-0.75)
\put(0.075,2.25){\vector(1,0){1.85}}
\put(2.375,1.925){\vector(0,-1){1.85}}
\put(-0.375,1.925){\vector(0,-1){1.85}}
\put(0.075,-0.25){\vector(1,0){1.85}}
\put(-0.375,2.25){\makebox(0,0){$X_\infty$}}
\put(2.375,2.25){\makebox(0,0){$X_\infty$}}
\put(-0.375,-0.25){\makebox(0,0){$X$}}
\put(2.375,-0.25){\makebox(0,0){$X$.}}
\put(1,2.375){\makebox(0,0)[b]{$\hat r$}}
\put(-0.25,1){\makebox(0,0)[l]{$\theta_n$}}
\put(1,-0.375){\makebox(0,0)[t]{$r$}}
\put(2.5,1){\makebox(0,0)[l]{$\theta_n$}}
\end{picture}
\end{equation*}
\par
Consider a probability measure $\mu$ on $X$ that satisfies
\begin{equation}\label{eqpr_4}
\int_Xf\,\mathrm{d}\mu=\int_X\frac{1}{\#r^{-1}(x)}\sum_{r(y)=x}f(y)\,\mathrm{d}\mu(x).
\end{equation}
\par
It is known that such measures $\mu$ on $X$ exist for a general
class of systems $r:X\rightarrow X$. The measure $\mu$ is said to
be {\it strongly $r$-invariant}. We have already discussed some in
Section \ref{Fun} above.
\par
If $X=X(A)$ is the state space of a sub-shift, we saw that
$\mu=\mu_A$ may be constructed as an application of Ruelle's
theorem (see Lemma \ref{lemfme3_2}). If $X=Julia(r)$ is the Julia
set of some rational mapping, then it is also known
\cite{Bea},\cite{Mil} that a strongly $r$-invariant measure $\mu$
on $X=Julia(r)$ exists.

\par

For $m_0\in L^\infty(X)$, define
\begin{equation}\label{eqpr_7}
(R\xi)(x)=\frac{1}{\#
r^{-1}(x)}\sum_{r(y)=x}|m_0(y)|^2\xi(y),\quad(\xi\in L^1(X)).
\end{equation}
\par
The next two theorems (Theorem \ref{thpr_1}-\ref{thpr_2}) are key
to our dilation theory. The dilations which we construct take
place at three levels as follows:
\begin{itemize}
\item Dynamical systems
\begin{equation*}(X,r,\mu)\mbox{ endomorphism }\rightarrow(X_\infty,\hat
r,\hat\mu),\mbox{ automorphism }.\end{equation*} \item Hilbert spaces
\begin{equation*}L_2(X,h\,\mathrm{d}\mu)\rightarrow (R_{m_0}h=h)\rightarrow
L^2(X_\infty,\hat\mu).\end{equation*} \item Operators
\begin{gather*}S_{m_0}\mbox{ isometry }\rightarrow U\mbox{ unitary (if
}m_0\mbox{ is non-singular});\\
M(g)\mbox{ multiplication operator }\rightarrow M_\infty(g).\end{gather*}
\end{itemize}
\begin{definition}
A function $m_0$ on a measure space is called {\it singular} if
$m_0=0$ on a set of positive measure.
\end{definition}
\par
In general, the operators $S_{m_0}$ on $H_0=L^2(X,h\,\mathrm{d}\mu)$, and
$U$ on $L^2(X_\infty,\hat\mu)$, may be given only by abstract
Hilbert space axioms; but in our {\it martingale representation},
we get the following two concrete formulas:
\begin{align*}\qquad(S_{m_0}\xi)(x)&=m_0(x)\xi(r(x)),& &(x\in X,\xi\in H_0);\\
(Uf)(\hat x)&=m_0(x_0)f(\hat r(\hat x)),& &(\hat x\in
X_\infty,f\in L^2(X_\infty,\hat\mu)).\qquad\end{align*}
\begin{theorem}\label{thpr_1}
If $h\in L^1(X)$, $h\geq0$ and $Rh=h$, then there exists a unique
measure $\hat\mu$ on $X_\infty$ such that
\begin{equation*}\hat\mu\circ\theta_n^{-1}=\omega_n,\quad(n\geq0),\end{equation*}
where \begin{equation}\label{eqpr_6}
\omega_n(f)=\int_XR^n(fh)\,\mathrm{d}\mu,\quad(f\in L^\infty(X)).
\end{equation}
\end{theorem}

\begin{proof}
It is enough to check that the measures $\omega_n$ and
$\omega_{n+1}$ are {\it compatible}, i.e., we have to check if
\begin{equation*}\omega_{n+1}(f\circ r)=\omega_n(f),\quad(f\in L^\infty(X)).\end{equation*}
But
\begin{equation*}R^{n+1}(f\circ r\,h)=R^n(R(f\circ r\,h))=R^n(fRh)=R^n(fh).\qedhere\end{equation*}
\end{proof}
Note that we can identify functions on $X$ with functions on
$X_\infty$ by
\begin{equation*}f(x_0,x_1...)=f(x_0),\quad(f:X\rightarrow\bc).\end{equation*}
\begin{theorem}\label{thpr_2}
\begin{equation}\label{eqpr_10}
\frac{d(\hat\mu\circ\hat r)}{d\hat\mu}=|m_0|^2
\end{equation}
\end{theorem}

\begin{proof}
Equation (\ref{eqpr_10}) can be rewritten as
\begin{equation*}\int_{X_\infty}|m_0|^2f\circ\hat
r\,\mathrm{d}\hat\mu=\int_{X_\infty}f\,\mathrm{d}\hat\mu,\quad(f\in
L^\infty(\hat\mu)).\end{equation*}
 By the uniqueness of $\hat\mu$, it is enough to check that
\begin{equation*}\int_{X_\infty}|m_0|^2(x_0)(f\circ\theta_n)\circ\hat
r(\hat x)\,\mathrm{d}\hat\mu(\hat x)=\omega_n(f),\quad(f\in L^\infty(X)),\end{equation*}
or, equivalently (since $\theta_n\hat r=r\theta_n$ and
$x_0=r^n(x_n)$):
\begin{equation}\label{eqtheorem53diez}\omega_n(|m_0|^2\circ r^n\,f\circ
r)=\omega_n(f).\end{equation} We can compute:
\begin{align*}&\int_XR^n(|m_0|^2\circ r^nf\circ r\,h)\,\mathrm{d}\mu=\int_X|m_0|^2R^n(f\circ r\,h)\,\mathrm{d}\mu\\
&\qquad=\int_X|m_0|^2R^{n-1}(fRh)\,\mathrm{d}\mu=\int_X|m_0|^2R^{n-1}(fh)\,\mathrm{d}\mu=\int_XR(R^{n-1}(fh))\,\mathrm{d}\mu,\end{align*}
and we used (\ref{eqpr_4}) for the last equality. This proves
(\ref{eqtheorem53diez}) and the theorem.
\end{proof}

\begin{theorem}\label{thpr_3}
Suppose $m_0$ is non-singular, i.e., it does not vanish on a set
of positive measure. Define $U$ on $L^2(X_\infty,\hat\mu)$ by
\begin{align*}Uf&=m_0f\circ\hat r,& &(f\in L^2(X_\infty,\hat\mu)),\\
\qquad\pi(g)f&=gf,& &(g\in L^\infty(X),f\in L^2(X_\infty,\hat\mu)),\qquad\\
\varphi&=1.& &\end{align*}
Then $(L^2(X_\infty,\hat\mu),U,\pi,\varphi)$ is the covariant
system associated to $m_0$ and $h$ as in Corollary \ref{cor2_4}.
Moreover, if $M_gf=gf$ for $g\in L^\infty(X_\infty,\hat\mu)$ and
$f\in L^2(X_\infty,\hat\mu)$, then
\begin{equation*}UM_gU^{-1}=M_{g\circ\hat r}.\end{equation*}
\end{theorem}

\begin{proof}
Theorem \ref{thpr_2} shows that $U$ is isometric. Since $m_0$ is
non-singular, the same theorem can be used to deduce that
\begin{equation*}U^*f=\frac{1}{m_0\circ\hat r^{-1}}f\circ\hat r^{-1}\end{equation*}
is a well defined inverse for $U$.
\par
The covariance relation follows by a direct computation. Also we
obtain
\begin{equation*}U^{-n}\pi(g)U^nf=g\circ\hat r^{-n}f,\quad(g\in L^\infty(X),f\in
L^2(X_\infty,\hat\mu)),\end{equation*} which shows that $\varphi$ is cyclic.
\par
The other requirements of Corollary \ref{cor2_4}, are easily
obtained by computation.
\end{proof}

\begin{remark}\label{rempr_4}
When $m_0$ is singular $U$ is just an isometry (not onto).
However, we still have many of the relations: the covariance
relation becomes
\begin{equation*}U\pi(f)=\pi(f\circ r)U,\quad(f\in L^\infty(X)),\end{equation*}
the scaling equation remains true,
\begin{equation}\label{eqpr4_1}
U\varphi=\pi(m_0)\varphi,
\end{equation}
and the correlation function of $\varphi$ is $h$:
\begin{equation*}\ip{\varphi}{\pi(f)\varphi}=\int_Xfh\,\mathrm{d}\mu,\quad(f\in
L^\infty(X)).\end{equation*} We further note that equation (\ref{eqpr4_1}) is
an abstract version of the scaling identity from wavelet theory.
In Section \ref{Intro} we recalled the scaling equation in its two
equivalent forms, the additive version (\ref{eqintsc1}), and its
multiplicative version (\ref{eqintsc2}). The two versions are
equivalent via the Fourier transform.
\end{remark}

\section{\label{Mar}\uppercase{Martingales}} We give now a different representation of
the construction of the covariant system associated to $m_0$ and
$h$ given in Theorem \ref{thpr_3}.
\par
Let
\begin{equation*}H_n:=\{f\in L^2(X_\infty,\hat\mu)\,|\,
f=\xi\circ\theta_n,\xi\in L^2(X,\omega_n)\}.\end{equation*} Then $H_n$ form an
increasing sequence of closed subspaces which have dense union.
\par
We can identify the functions in $H_n$ with functions in
$L^2(X,\omega_n)$, by
\begin{equation*}i_n(\xi)=\xi\circ\theta_n,\quad(\xi\in L^2(X,\omega_n)).\end{equation*}
The definition of $\hat\mu$ makes $i_n$ an isomorphism between
$H_n$ and $L^2(X,\omega_n).$
\par
Define
\begin{equation*}\mathcal{H}:=\{(\xi_0,\xi_1,...)\,|\, \xi_n\in
L^2(X,\omega_n),R(\xi_{n+1}h)=\xi_nh,\,
\sup_n\int_XR^n(|\xi_n|^2h)\,\mathrm{d}\mu<\infty\},\end{equation*} with the scalar
product
\begin{equation*}\ip{(\xi_0,\xi_1,...)}{(\eta_0,\eta_1,...)}=\lim_{n\rightarrow\infty}\int_XR^n(\cj\xi_n\eta_nh)\,\mathrm{d}\mu.\end{equation*}
\begin{theorem}\label{thpr_5}
The map $\Phi:L^2(X_\infty,\hat\mu)\rightarrow\mathcal{H}$ defined
by
\begin{equation*}\Phi(f)=(i_n^{-1}(P_nf))_{n\geq0},\end{equation*}
where $P_n$ is the projection onto $H_n$ is an isomorphism.
\begin{align*}\Phi U\Phi^{-1}(\xi_n)_{n\geq_0}&=(m_0\circ
r^n\,\xi_{n+1})_{n\geq0},\\
\Phi\pi(g)\Phi^{-1}(\xi_n)_{n\geq_0}&=(g\circ
r^n\,\xi_n)_{n\geq0},\\
\Phi\varphi&=(1,1,...).\vphantom{\Phi\pi(g)\Phi^{-1}(\xi_n)_{n\geq_0}}\end{align*}
\end{theorem}
\begin{proof}
Let $\xi_n:=i_n^{-1}(P_nf).$ We check that $R(\xi_{n+1}h)=\xi_nh$.
For this it is enough to see that the projection of
$\xi_{n+1}\circ\theta_{n+1}$ onto $H_n$ is
$(R(\xi_{n+1}h)/h)\circ\theta_n$. We compute the scalar products
with $g\circ\theta_n\in H_n$:
\begin{align*}\ip{\xi_{n+1}\circ\theta_{n+1}}{g\circ\theta_n}&=\int_{X_\infty}
\!
\cj\xi_{n+1}\circ\theta_{n+1}g\circ
r\circ\theta_{n+1}\,\mathrm{d}\hat\mu=\int_X
\!
R^{n+1}(\cj\xi_{n+1}g\circ
rh)\,\mathrm{d}\mu\\
&=\int_X
\!
R^n(g\frac{R(\cj\xi_{n+1}h)}{h}h)\,\mathrm{d}\mu=\Bip{\frac{R(\xi_{n+1}h)}{h}\circ\theta_n}{g\circ\theta_n}.\end{align*}
Since the union of $(H_n)$ is dense, $P_nf$ converges to $f$. As
each $i_n$ is isometric,
\begin{equation*}\ip{f}{g}=\lim_{n\rightarrow\infty}\ip{P_nf}{P_ng}=\lim_{n\rightarrow\infty}\ip{\Phi(f)_n}{\Phi(g)_n}_{L^2(X,\omega_n)}=\ip{\Phi(f)}{\Phi(g)}.\end{equation*}
\par
Now we check that $\Phi$ is onto. Take
$(\xi_n)_{n\geq0}\in\mathcal{H}$. Then define
\begin{equation*}f_n:=\xi_n\circ\theta_n=i_n^{-1}(\xi_n).\end{equation*}
The previous computation shows that \begin{equation*}P_{n}f_{n+1}=f_n.\end{equation*} Also
\begin{equation*}\sup_n\|f_n\|^2=\sup_n\int_XR^n(|\xi_n|^2h)\,\mathrm{d}\mu<\infty.\end{equation*}
But then, by a standard Hilbert space argument, $f_n$ is a Cauchy
sequence which converges to some
\begin{equation*}f=\lim_{n\rightarrow\infty}f_n=f_0+\sum_{k=0}^\infty(f_{k+1}-f_k)\in L^2(X_\infty,\mu)\end{equation*} with
$P_nf=f_n$ for all $n\geq 0$, and we conclude that
$\Phi(f)=(\xi_n)_{n\geq0}$.
\par
The form of  $\Phi U\Phi^{-1}$ and $\Phi\pi(g)\Phi^{-1}$ can be
obtain from the next lemma (using the fact that $P_nUf=UP_{n+1}$).
\end{proof}
\begin{lemma}\label{lemmart6}

The following diagram is commutative
\begin{equation*}
\begin{array}{ccc}
L^2(X,\omega_n) & \overset{\displaystyle \alpha}{\longrightarrow} & L^2(X,\omega_{n+1}) \\[3\jot]
\downarrow\rlap{$\mkern2mu i_n$} & & \downarrow\rlap{$\mkern2mu i_{n+1}$} \\[3\jot]
H_n & \hooklongrightarrow & H_{n+1}
\end{array}
\end{equation*}
where $\alpha(\xi)=\xi\circ r$.
\par
If $\xi\circ\theta_{n+k}\in H_{n+k}$, then
\begin{equation}\label{eqmart4_2}
P_n(\xi\circ\theta_{n+k})=\frac{R^k(\xi h)}{h}\circ\theta_n.
\end{equation}
\begin{equation}\label{eqmart4_4}
U^*f=\chi_{\{m_0\circ\hat r^{-1}\neq0\}}\frac{1}{m_0\circ\hat
r^{-1}}f\circ\hat r^{-1},\quad(f\in L^2(X_\infty,\hat\mu)).
\end{equation}

\begin{equation}\label{eqmart4_3}
UP_{n+1}U^*=P_n,\quad(n\geq 0).
\end{equation}
\end{lemma}

\begin{proof}
For $\xi\in L^2(X,\omega_n)$, $\xi\circ\theta_n=\xi\circ
r\circ\theta_{n+1}=i_{n+1}(\alpha(\xi))$, thus the diagram
commutes.
\par
We have to check that, for all $\eta\in L^2(X,\omega_n)$ we have
\begin{equation*}\ip{\xi\circ\theta_{n+k}}{\eta\circ\theta_n}=\Bip{\frac{R^k(\xi
h)}{h}\circ\theta_n}{\eta\circ\theta_n}.\end{equation*} But
\begin{align*}\ip{\xi\circ\theta_{n+k}}{\eta\circ\theta_n}&=\int_XR^{n+k}(\cj\xi\eta\circ
r^k\,h)\,\mathrm{d}\mu=\int_XR^n(\frac{R^k(\cj\xi h)}{h}\eta h)\,\mathrm{d}\mu\\
&=\Bip{\frac{R^k(\xi h)}{h}\circ\theta_n}{\eta\circ\theta_n}.\end{align*}
\par
Equation (\ref{eqmart4_4}) can be proved by a direct computation.
\par
Since $(H_n)$ are dense in $L^2(X_\infty,\hat\mu)$, we can check
(\ref{eqmart4_3}) on $H_{n+k}$. Take $\xi\circ\theta_{n+k}\in
H_{n+k}$, then
\begin{align*}
UP_{n+1}U^*(\xi\circ\theta_{n+k})&=UP_{n+1}\left(\chi_{\{m_0\circ\hat
r^{-1}\neq0\}}\frac{1}{m_0\circ\hat
r^{-1}}\xi\circ\theta_{n+k}\circ\hat r^{-1}\right)\\
=UP_{n+1}&\left(\left(\chi_{\{m_0\circ\hat r^{-1}\neq0\}}\circ
r^{n+k+1}\frac{1}{m_0\circ
r^{n+k}}\xi\right)\circ\theta_{n+k+1}\right)\\
&=U\left(\left(\frac{R^{k}\left(\chi_{\{m_0\circ\hat
r^{-1}\neq0\}}\circ r^{n+k+1}\frac{1}{m_0\circ r^{n+k}}\xi
h\right)}{h}\right)\circ\theta_{n+1}\right)\\
&=U\left(\left(\chi_{\{m_0\circ\hat r^{-1}\neq0\}}\circ
r^{n+1}\frac{1}{m_0\circ r^n}\frac{R^{k}(\xi
h)}{h}\right)\circ\theta_{n+1}\right)\\
&=m_0\chi_{\{m_0\circ\hat r^{-1}\neq0\}}\circ
r\frac{1}{m_0}\frac{R^{k}(\xi h)}{h}\circ\theta_n\\
&=P_n(\xi\circ\theta_{n+k}).\qedhere
\end{align*}
\end{proof}
As a consequence of Lemma \ref{lemmart6} we also have:
\begin{proposition}\label{propmart7}
The identification of functions in $L^2(X,\omega_n)$ with
martingales is given by
\begin{equation}\label{eqmart7_1}
\Phi(i_n(\xi))=\left(\frac{R^n(\xi h)}{h},...,\frac{R(\xi
h)}{h},\xi,\xi\circ r,\xi\circ r^2,...\right),\quad(\xi\in
L^2(X,\omega_n),n\geq0).
\end{equation}
\end{proposition}
\par
The condition that $m_0$ be non-singular is essential if one wants
$U$ to be unitary. We illustrate this by an example.
\begin{example}[Shannon's wavelet]
\par
Let $\br/\bz\simeq[-\frac12,\frac12)$. By this we mean that
functions on  $[-\frac12,\frac12)$ are viewed also as functions on
$\br$ via periodic extension, i.e., $f(x+n)=f(x)$ if
$x\in[-\frac12,\frac12)$ and $n\in\bz$.
\par
Set
\begin{equation*}m_0(x)=\sqrt{2}\chi_{[-\frac14,\frac14)}(x).\end{equation*}
Then \begin{equation}\label{eqmartex8_1}
\hat\varphi(x)=\prod_{k=1}^\infty\frac{1}{\sqrt{2}}m_0\left(\frac{x}{2^k}\right)=\chi_{[-\frac14,\frac14)}\left(\frac{x}2\right)=\chi_{[-\frac12,\frac12)}(x),
\end{equation}
and
\begin{equation*}\varphi(t)=\frac{\sin\pi t}{\pi t}.\end{equation*}
For functions in $L^1(\br/\bz)$, the Ruelle operator $R_{m_0}$ is
\begin{align*}(R_{m_0}f)(x)&=\chi_{[-\frac14,\frac14)}(\frac{x}{2})f(\frac{x}{2})+\chi_{[-\frac14,\frac14)}(\frac{x+1}{2})f(\frac{x+1}{2})=\chi_{[-\frac14,\frac14)}(\frac{x}{2})f(\frac{x}{2})\\
&=\chi_{[-\frac12,\frac12)}(x)f(\frac{x}{2})=f(\frac{x}{2}),\mbox{
for }x\in[-\frac12,\frac12).\end{align*} Hence $R_{m_0}1=1$.
\par
Note from (\ref{eqmartex8_1}) that $\hat\varphi(x+n)=0$ if
$n\in\bz\setminus\{0\}$.\par Let $\xi\in L^2(\br/\bz)$. Then we
get
\begin{align*}\int_{X_\infty}{|\xi\circ\theta_n|^2}\,\mathrm{d}\hat\mu&=\int_X|\xi|^2\,\mathrm{d}\omega_n=\int_{-\frac12}^{\frac12}R^n(|\xi|^2)(x)\,\mathrm{d}x\\
&=\int_{-\frac12}^{\frac12}|\xi(2^{-n}x)|^2\,\mathrm{d}x=2^n\int_{-\frac{1}{2^{n+1}}}^{\frac{1}{2^{n+1}}}|\xi(x)|^2\,\mathrm{d}x.\end{align*}
But then
$L^2(X,\omega_n)=L^2([-\frac{1}{2^{n+1}},\frac{1}{2^{n+1}}),2^n\,\mathrm{d}x)$
and we see that the map
 \begin{equation*}\alpha:L^2(X,\omega_n)\rightarrow
L^2(X,\omega_{n+1}),\quad \alpha(\xi)=\xi(2\cdot)\end{equation*} is an isometry
(Lemma \ref{lemmart6}) which is also surjective with inverse
$\xi\mapsto\xi(\frac{x}{2})$.
\par
With Lemma \ref{lemmart6}, we get that the inclusion of $H_n$ in
$H_{n+1}$ is in fact an identity, therefore
\begin{equation*}L^2(X_\infty,\hat\mu)=H_0=L^2([-\frac12,\frac12),\,\mathrm{d}x).\end{equation*}
\end{example}

\par
When $m_0$ is non-singular, Theorem \ref{thpr_3} shows that the
covariant system\\ $(L^2(X_\infty,\hat\mu),U,\pi,\varphi)$ has $U$
unitary so, by uniqueness, it is isomorphic to the one constructed
via the Kolmogorov theorem in Corollary \ref{cor2_4}, which we
denote by $(\tilde H,\tilde U,\tilde\pi,\tilde\varphi)$.
\par
The next theorem shows that even when $m_0$ is singular, the
covariant system
\linebreak
$(L^2(X_\infty,\hat\mu),U,\pi,\varphi)$ can be
embedded in the $(\tilde H,\tilde U,\tilde\pi,\tilde\varphi)$.
\begin{theorem}\label{thcond3}
There exists a unique isometry
$\Psi:L^2(X_\infty,\hat\mu)\rightarrow\tilde H$ such that
\begin{equation*}\Psi(\xi\circ\theta_n)=\tilde U^{-n}\tilde\pi(\xi)\tilde
U^n\tilde\varphi,\quad(\xi\in L^\infty(X,\mu)).\end{equation*} $\Psi$
intertwines the two systems, i.e.,
\begin{equation*}\Psi U=\tilde U\Psi,\quad \Psi\pi(g)=\tilde\pi(g)\Psi,\mbox{ for }g\in
L^\infty(X,\mu),\quad\Psi\varphi=\tilde\varphi.\end{equation*}
\end{theorem}

\begin{proof}
Let $j_n:H_n\rightarrow\tilde H$ be defined on a dense subspace by
\begin{equation*}j_n(\xi\circ\theta_n)=\tilde U^{-n}\tilde\pi(\xi)\tilde
U^n\tilde\varphi,\quad(\xi\in L^\infty(X,\mu)).\end{equation*} Then $j_n$ is a
well defined isometry because
\begin{equation*}\|\xi\circ\theta_n\|_{L^2(\hat\mu)}^2=\int_XR^n(|\xi|^2h)\,\mathrm{d}\mu=\int_X|m_0^{(n)}|^2|\xi|^2\,\mathrm{d}\mu=\|\tilde U^{-n}\tilde\pi(\xi)\tilde
U^n\tilde\varphi\|^2,\end{equation*} where \begin{equation*}m_0^{(n)}:=m_0\cdot m_0\circ
r\cdot..\cdot m_0\circ r^{n-1}.\end{equation*} Also note that
\begin{equation*}j_{n+1}(\xi\circ\theta_n)=j_n(\xi\circ r\circ\theta_{n+1})=\tilde U^{-n-1}\tilde\pi(\xi\circ r)\tilde
U^{n+1}\tilde\varphi=\tilde U^{-n}\tilde\pi(\xi)\tilde
U^n\tilde\varphi,\end{equation*} so we can construct $\Psi$ on
$L^2(X_\infty,\hat\mu)$ such that it agrees with $j_n$ on $H_n$.
\par
Next, we check the intertwining properties; it is enough to verify
them on $H_n$:
\begin{align*}
\tilde U\Psi(\xi\circ\theta_n)&=\tilde U\tilde
U^{-n}\tilde\pi(\xi)\tilde U^n\tilde\varphi=\tilde
U^{-n+1}\tilde\pi(\xi\circ r)\tilde U^{n-1}\tilde
U\tilde\varphi\\
&=\tilde U^{-n+1}\tilde\pi(\xi)\tilde
U^{n-1}\tilde\pi(m_0)\tilde\varphi,\\
\Psi U(\xi\circ\theta_n)&=\Psi(m_0\xi\circ\theta_n\circ\hat
r)=\Psi((m_0\circ r^{n-1}\xi)\circ\theta_{n-1})\\
&=\tilde
U^{-n+1}\tilde\pi(m_0\circ r^{n-1}\xi)\tilde
U^{n-1}\tilde\varphi=\tilde U^{-n+1}\tilde\pi(\xi)\tilde
U^{n-1}\tilde\pi(m_0)\tilde\varphi.\end{align*} The other intertwining
relations can be checked by some similar computations.
\end{proof}

\subsection{\label{Con}Conditional expectations}
We can consider the $\sigma$-algebras
\begin{equation*}\mathfrak{B}_n:=\theta_n^{-1}(\mathfrak{B}),\end{equation*} $\mathfrak{B}$
being the $\sigma$-algebra of Borel subsets in $X$. Note that
$\theta_n^{-1}(E)=\theta_{n+1}^{-1}(r^{-1}(E))$. If follows that
\begin{equation*}\mathfrak{B}_0\subset\mathfrak{B}_1\subset...\subset\mathfrak{B}_n\subset\mathfrak{B}_{n+1}\subset...\end{equation*}
We set $\mathfrak{B}_\infty=\cup_{n\geq0}\mathfrak{B}_n$ which is
a sigma-algebra on $X_\infty$.
\par
The functions on $X_\infty$ which are $\mathfrak{B}_n$ measurable
are the functions which depend only on $x_0,...,x_n$. $H_n$
consists of function in $L^2(X_\infty,\mathfrak{B}_n,\hat\mu).$
Also we can regard $L^\infty(X_\infty,\mathfrak{B}_n,\hat\mu)$ as
an increasing sequence of subalgebras of
$L^\infty(X_\infty,\hat\mu)$. The map
\begin{equation*}i_n: L^\infty(X,\omega_n)\rightarrow
L^\infty(X_\infty,\mathfrak{B}_n,\hat\mu)\end{equation*} is an isomorphism.
\par
An application of the Radon-Nikodym theorem shows that there
exists a unique {\it conditional expectation}
$E_n:L^1(X_\infty,\hat\mu)\rightarrow
L^1(X_\infty,\mathfrak{B}_n,\hat\mu)$ determined by the relation
\begin{equation}\label{eqcond1}
\int_{X_\infty}E_n(f)g\,\mathrm{d}\hat\mu=\int_{X_\infty}fg\,\mathrm{d}\hat\mu,\quad(g\in
L^\infty(X_\infty,\mathfrak{B}_n,\hat\mu)).
\end{equation}
We enumerate the properties of these conditional expectations.
\begin{proposition}\label{propcond2}
\begin{equation}\label{eqcond2_1}
E_n(fg)=fE_n(g),\quad(f\in
L^\infty(X_\infty,\mathfrak{B}_n,\hat\mu),g\in
L^1(X_\infty,\hat\mu)),
\end{equation}
\begin{equation}\label{eqcond2_2}
E_n(f)\geq 0,\mbox{ if }f\geq 0,
\end{equation}
\begin{equation}\label{eqcond2_3}
E_mE_n=E_nE_m=E_n,\mbox{ if }m\geq n,
\end{equation}
\begin{equation}\label{eqcond2_4}
\int_{X_\infty}E_n(f)\,\mathrm{d}\hat\mu=\int_{X_\infty}f\,\mathrm{d}\hat\mu,
\end{equation}
\begin{equation}\label{eqcond2_5}
E_n(f)=P_n(f),\mbox{ if }f\in L^2(X_\infty,\hat\mu).
\end{equation}
\end{proposition}
\begin{definition}
A sequence $(f_n)_{n\geq0}$ of measurable functions on $X_\infty$
is said to be a {\it martingale} if
\begin{equation*}E_nf_{n+1}=f_n,\quad(n\geq0),\end{equation*}
where $E_n$ is a family of conditional expectations as in
Proposition \ref{propcond2}.
\end{definition}

\begin{proposition}\label{propcond3}
If $\xi\in L^1(X,\omega_{n+k})$ then
\begin{equation}\label{eqcond3_1}
E_n(\xi\circ\theta_{n+k})=\frac{R^k(\xi h)}{h}\circ\theta_n
\end{equation}
\end{proposition}
\begin{proof}
If $\xi\in L^2(X,\omega_n)$, the formula follows from Lemma
\ref{lemmart6}. The rest follows by approximation.
\end{proof}
\par
Proposition \ref{propcond3} offers a direct link between the
operator powers $R^k$ and the conditional expectations $E_n$. It
shows in particular how our martingale construction depends on the
Ruelle operator $R$. For a sequence $(\xi_n)_{n\geq0}$ of
measurable functions on $X$, $(\xi_n\circ\theta_n)_{n\geq0}$ is a
martingale if and only if
\begin{equation*}R(\xi_{n+1}h)=\xi_nh,\quad(n\geq0).\end{equation*}
\par
A direct application of Doob's theorem (Theorem IV-1-2, in
\cite{Neveu}) gives the following:
\begin{proposition}\label{propcond4}
If $\xi_n\in L^1(X,\omega_n)$ is a sequence of functions with the
property that
\begin{equation*}R(\xi_{n+1}h)=\xi_nh,\quad(n\geq0),\end{equation*}
then the sequence $\xi_n\circ\theta_n$ converges $\hat\mu$-almost
everywhere.
\end{proposition}
Then Proposition IV-2-3 from \cite{Neveu}, translates into
\begin{proposition}\label{propcond5}
Suppose $\xi_n\in L^1(X,\omega_n)$ is a sequence with the property
that
\begin{equation*}R(\xi_{n+1}h)=\xi_nh,\quad(n\geq0).\end{equation*}
The following conditions are equivalent:
\begin{enumerate}
\item The sequence $\xi_n\circ\theta_n$ converges in
$L^1(X_\infty,\hat\mu)$. \item
$\sup_n\int_XR^n(|\xi|h)\,\mathrm{d}\mu<\infty$ and the a.e. limit
$\xi_\infty=\lim_n\xi_n\circ\theta_n$ satisfies
$\xi_n\circ\theta_n=E_n(\xi_\infty).$ \item There exists a
function $\xi\in L^1(X_\infty,\hat\mu)$ such that
$\xi_n\circ\theta_n=E_n(\xi)$ for all $n$. \item The sequence
$\xi_n\circ\theta_n$ satisfies the uniform integrability
condition:
\begin{equation*}\sup_n\int_XR^n(\chi_{\{|\xi_n|>a\}}\xi_nh)\,\mathrm{d}\mu\downarrow0\mbox{
as }a\uparrow\infty\end{equation*}.
\end{enumerate}
If one of the conditions is satisfied, the martingale $(\xi_n)_n$
is called regular.
\end{proposition}
\par
Convergence in $L^p$ is given by Proposition IV-2-7 in
\cite{Neveu}:
\begin{proposition}\label{propcond6} Let $p>1$. Every martingale
$(\xi_n)_n$ with $\xi_n\in L^p(X,\omega_n)$ and
\begin{equation*}\sup_n\|\xi_n\|_p<\infty\end{equation*}
is regular, and $\xi_n\circ\theta_n$ converges in
$L^p(X_\infty,\hat\mu)$ to $\xi_\infty$.
\end{proposition}
\par
We have seen that functions $f$ on $X_\infty$ may be identified
with sequences $(\xi_n)$ of functions on $X$. When $r:X\rightarrow
X$ is given, the induced mappings
\begin{equation}\label{eqmartins1}
\hat r:X_\infty\rightarrow X_\infty,\mbox{ and }\hat
r^{-1}:X_\infty\rightarrow X_\infty
\end{equation}
yield transformations of functions on $X_\infty$ as follows
$f\mapsto f\circ\hat r$ and $f\mapsto f\circ\hat r^{-1}$.
\par
The 1-1 correspondence
\begin{equation}\label{eqmartins2}
f\mbox{ function on }X_\infty\leftrightarrow\xi_0,\xi_1,...\mbox{
functions on }X
\end{equation}
is determined uniquely by
\begin{equation}\label{eqmart3}
E_n(f)=\xi_n\circ\theta_n,\quad n=0,1,...
\end{equation}
When $f$ and $h$ are given, then the functions $(\xi_n)$ in
(\ref{eqmartins2}) must satisfy
\begin{equation}\label{eqmartins4}
R(\xi_{n+1}h)=\xi_nh,\quad(n\geq0)
\end{equation}
\begin{proposition}\label{propmartins4_11}
Assume $m_0$ is non-singular. If $f$ is a function on $X_\infty$
and $f\leftrightarrow(\xi_n)$ as in \textup{(\ref{eqmartins2})} then
\begin{equation}\label{eqmartins5}
f\circ\hat r\leftrightarrow\xi_{n+1}
\end{equation}
\begin{equation}\label{eqmartins6}
f\circ\hat r^{-1}\leftrightarrow\xi_{n-1}
\end{equation}
Specifically we have
\begin{equation}\label{eqmartins7}
E_n(f\circ\hat r)=\xi_{n+1}\circ\theta_n
\end{equation}
and
\begin{equation}\label{eqmartins8}
E_n(f\circ\hat
r^{-1})=\xi_{n-1}\circ\theta_n=\left(\frac{R(\xi_nh)}{h}\right)\circ\theta_n
\end{equation}
Or equivalently
\begin{equation}\label{eqmartins9}
f\circ\hat r\leftrightarrow(\xi_1,\xi_2,...),
\end{equation}
and
\begin{equation}\label{eqmartins10}
f\circ\hat
r^{-1}\leftrightarrow(\frac{R(\xi_0h)}{h},\xi_0,\xi_1,...).
\end{equation}
\end{proposition}
\begin{proof}Theorem \ref{thpr_2} is used in both parts of the
proof below. 

We have for $g:X\rightarrow\bc$,
\begin{align*}
\int_{X_\infty}\!E_n(f\circ\hat
r)\,g\circ\theta_n\,\mathrm{d}\hat\mu&=\int_{X_\infty}\!f\circ\hat
r\,g\circ\theta_{n+1}\circ\hat r\,\mathrm{d}\hat\mu
=\int_{X_\infty}\!\frac{1}{|m_0|^2\circ\hat
r^{-1}}f\,g\circ\theta_{n+1}\,\mathrm{d}\hat\mu\\
&=\int_{X_\infty}\!E_{n+1}(f)\left(\frac{1}{|m_0|^2\circ
r^n}g\right)\circ\theta_{n+1}\,\mathrm{d}\hat\mu\\
&=\int_{X_\infty}\!\xi_{n+1}\circ\theta_{n}\circ\hat
r^{-1}\left(\frac{1}{|m_0|^2\circ
r^n}g\right)\circ\theta_{n}\circ\hat r^{-1}\,\mathrm{d}\hat\mu\\
&=\int_{X_\infty}\!|m_0|^2\xi_{n+1}\circ\theta_n\frac{1}{|m_0|^2}\,g\circ\theta_n\,\mathrm{d}\hat\mu\\
&=\int_{X_\infty}\!\xi_{n+1}\circ\theta_n\,g\circ\theta_n\,\mathrm{d}\hat\mu.
\end{align*}
Thus $E_n(f\circ\hat r)=\xi_{n+1}\circ\theta_n$.
\par
\begin{align*}
\int_{X_\infty}\!E_n(f\circ\hat
r^{-1})g_n\circ\theta_n\,\mathrm{d}\hat\mu&=\int_{X_\infty}\!f\circ\hat
r^{-1}\,g_n\circ\theta_{n-1}\circ\hat r^{-1}\,\mathrm{d}\hat\mu\\
&=\int_{X_\infty}\!|m_0|^2f\,g_n\circ\theta_{n-1}\,\mathrm{d}\hat\mu\\
&=\int_{X_\infty}\!E_{n-1}(f)\left(|m_0|^2\circ
r^{n-1}\,g\right)\circ\theta_{n-1}\,\mathrm{d}\hat\mu\\
&=\int_{X_\infty}\!\xi_{n-1}\circ\theta_{n}\circ\hat
r\left(|m_0|^2\circ
r^{n-1}\,g\right)\circ\theta_{n}\circ\hat r\,\mathrm{d}\hat\mu\\
&=\int_{X_\infty}\!\frac{1}{|m_0|^2\circ\hat
r^{-1}}\xi_{n-1}\circ\theta_n\,|m_0|^2\circ\hat
r^{-1}\,g\circ\theta_n\,\mathrm{d}\hat\mu\\
&=\int_{X_\infty}\!\xi_{n-1}\circ\theta_n\,g\circ\theta_n\,\mathrm{d}\hat\mu
\end{align*}
and this implies (\ref{eqmartins8}).
\end{proof}

\section{\label{Int}\uppercase{Intertwining operators and cocycles}}
\par
 In the paper \cite{DaLa98}, Dai and Larson showed that the familiar
orthogonal wavelet systems have an attractive representation
theoretic formulation. This formulation brings out the geometric
properties of wavelet analysis especially nicely, and it led to
the discovery of {\it wavelet sets}, i.e., singly generated
wavelets in $L^2(\br^d)$, i.e., $\psi\in L^2(\br^d)$ such that
$\hat\psi=\chi_E$ for some $E\subset\br^d$, and
\begin{equation*}\{|\mbox{det}A|^{j/2}\psi(A^j\cdot-k)\,|\,j\in\bz,k\in\bz^d\}\end{equation*}
is an orthonormal basis.
\par
 The case when the initial resolution subspace for some wavelet
construction is singly generated, the wavelet functions should be
thought of as {\it wandering vectors}. If the scaling operation is
realized as a unitary operator $U$ in the Hilbert space
$H:=L^2(\br^d)$, then the notion of wandering, refers to vectors,
or subspaces which are mapped into orthogonal vectors (
respectively, subspaces) under powers of $U$. Since this approach
yields wavelet bases derived directly from the initial data, i.e.,
from the wandering vectors, $U$, and the integral translations,
the question of intertwining operators is a natural one. The
initial data defines a representation $\rho$.
\par
 An operator in $H$ which intertwines $\rho$ with itself is said to be in the
commutant of $\rho$; and Dai and Larson gave a formula for the
commutant. They showed that the operators in the commutant are
defined in a natural way from a class of invariant bounded
measurable functions, called {\it wavelet multipliers}. This and
other related results can be shown to generalize to the case of
operators which intertwine two wavelet representations $\rho$ and
$\rho'$.
\par
 Since our present martingale construction is a  generalization of the
traditional wavelet resolutions, see \cite{Jor06}, it is natural to ask for
theorems which generalize the known theorems about wavelet functions. We
prove in this section such a theorem, Theorem \ref{thintertw}. The applications of this
are manifold, and include the projective systems defined from Julia sets,
and from the state space of a subshift in symbolic dynamics.
\par
    Our formula for the commutant in this general context of projective
systems is shown to be related to the Perron-Frobenius-Ruelle
operator in Corollary \ref{corinter2}. This result implies in
particular that the commutant is abelian; and it makes precise the
way in which the representation $\rho$ itself decomposes as a
direct integral over the commutant. \par
Our proof of this
corollary depends again on Doob's martingale convergence theorem,
see (\ref{eqinterc3}) below, Section \ref{Mar} above, and \cite{Jor06},
Chapter 2.7.

\begin{definition}
If $m_0\in L^\infty(X)$ and $h\in L^1(X)$, we call $(m_0,h)$ a
{\it Perron-Ruelle-Frobenius pair }if
\begin{equation*}R_{m_0}h=h.\end{equation*}
\end{definition}
\begin{theorem}\label{thintertw}
Let $(m_0,h)$ and $(m_0',h')$ be two Perron-Ruelle-Frobenius pairs
with $m_0,m_0'$ non-singular, and let
$(L^2(X_\infty,\hat\mu),U,\pi,\varphi)$,
$(L^2(X_\infty,\hat\mu'),U',\pi',\varphi')$ be the associated
covariant systems. Let $X_\infty=X_\infty^a\cup X_\infty^s$ be the
Jordan decomposition of $\hat\mu'$ with respect to $\hat\mu$,
$X_\infty^a\cap X_\infty^s=\emptyset$, with
$\hat\mu(X_\infty^s)=0$ and $\hat\mu'|_{X_\infty^a}\prec\hat\mu$,
and denote by
\begin{equation*}\Delta:=\frac{d\,\hat\mu'|_{X_\infty^a}}{d\hat\mu}.\end{equation*}
Then there is a 1-1 correspondence between each two of the
following sets of data:
\begin{enumerate}
\item Operators $A:L^2(X_\infty,\hat\mu)\rightarrow
L^2(X_\infty,\hat\mu')$ that intertwine the covariant system,
i.e.,
\begin{equation}\label{eqinter1}
U'A=AU,\mbox{ and }\pi'(g)A=\pi(g)A,\mbox{ for }g\in L^\infty(X).
\end{equation}
\item $\mathfrak{B}_\infty$-measurable functions
$f:X\rightarrow\bc$ such that $f|_{X_\infty^s}=0$,
$f\Delta^{\frac{1}{2}}$ is $\hat\mu$-bounded and
\begin{equation}\label{eqinter2}
m_0f=m_0'f\circ\hat r,\quad\hat\mu'-\mbox{a.e.}
\end{equation}
\item Measurable functions $h_0:X\rightarrow\bc$ such that
\begin{equation}\label{eqinter3_0}
|h_0|^2\leq chh'\,\mu\mbox{-a.e.},
\end{equation}
for some finite constant $c\geq 0$, with
\begin{equation}\label{eqinter3}
\frac{1}{\#r^{-1}(x)}\sum_{r(y)=x}\cj{m_0'(y)}m_0(y)h_0(y)=h_0(x),\mbox{
for }\mu\mbox{-a.e. }x\in X.
\end{equation}
{}From \textup{(i)} to \textup{(ii)} the correspondence is given by
\begin{equation}\label{eqinter4}
A\xi=f\xi,\quad(\xi\in L^2(X_\infty,\hat\mu)).
\end{equation}
{}From \textup{(ii)} to \textup{(iii)}, the correspondence is given by
\begin{equation}\label{eqinter5}
h_0=E_0^{\hat\mu'}(f)h'=E_0^{\hat\mu}(f\Delta)h
\end{equation}
{}From \textup{(i)} to \textup{(iii)} the correspondence is given by
\begin{equation}\label{eqinter6}
\ip{\varphi'}{A\pi(g)\varphi}=\int_Xgh_0\,\mathrm{d}\mu,\quad(g\in
L^\infty(X)).
\end{equation}
\end{enumerate}
\end{theorem}

\begin{proof}
Take $A$ as in (i). Then for all $g\in L^\infty(X)$ and any
$n\geq0$ we have that
\begin{equation*}A(g\circ\hat r^{-n})=A(U^{-n}\pi(g)U^n)(1)=(U'^{-n}\pi'(g)U'^n)(A(1))=g\circ\hat r^{-n}\cdot(A(1)).\end{equation*}
Denote by $f:=A(1)\in L^2(X_\infty,\hat\mu').$
\par
Since any $\mathfrak{B}_\infty$-measurable, bounded function
$\xi:X_\infty\rightarrow\bc$ can be pointwise $\hat\mu-$ and
$\hat\mu'-$approximated by functions of the form $g\circ\hat
r^{-n}$, we get that
\begin{equation*}A(\xi)=f\xi.\end{equation*}
We have also that
\begin{equation*}\int_{X_\infty}|f|^2|\xi|^2\,\mathrm{d}\hat\mu'\leq\|A\|^2\int_{X_\infty}|\xi|^2\,\mathrm{d}\hat\mu\end{equation*}
so
\begin{equation*}\int_{X_\infty^a}|f|^2|\xi|^2\Delta\,\mathrm{d}\hat\mu+\int_{X_\infty^s}|f|^2|\xi|^2\,\mathrm{d}\hat\mu'\leq\|A\|^2\int_{X_\infty}|\xi|^2\,\mathrm{d}\hat\mu\end{equation*}
Taking $\xi=\chi_{X_\infty^s}$ we obtain that $f=0$
$\hat\mu'$-a.e. on $X_\infty^s$; so we may take $f=0$ on
$X_\infty^s$. Then we get also that $|f\Delta^{1/2}|\leq\|A\|$
$\hat\mu$-a.e.
\par
Then, again by approximation we obtain that
\begin{equation*}A\xi=f\xi,\mbox{ for }\xi\in L^2(X_\infty,\hat\mu).\end{equation*}
\par
We have in addition the fact that $U'A=AU$, and this implies
(\ref{eqinter2}).
\par
Conversely, the previous calculations show that any operator
defined by (\ref{eqinter4}) with $f$ as in (ii), will be a bounded
operator which intertwines the covariant systems.
\par
Now take $A$ as in (i) and consider the linear functional
\begin{equation*}g\in L^\infty(X)\mapsto\ip{\varphi'}{A\pi(g)\varphi}\end{equation*}
This defines a measure on $X$ which is absolutely continuous with
respect to $\mu$. Let $h_0$ be its Radon-Nikodym derivative. We
have
\begin{align*}
\int_Xgh_0\,\mathrm{d}\mu&=\ip{\varphi'}{A\pi(g)\varphi}=
\ip{U'\varphi'}{U'A\pi(g)\varphi}\\
&=\ip{\pi'(m_0')\varphi'}{A\pi(g\circ r)\pi(m_0)\varphi}=\int_X\cj{m_0'}m_0g\circ r\,h_0\,\mathrm{d}\mu\\
&=\int_Xg\frac{1}{\#r^{-1}(x)}\sum_{r(y)=x}\cj{m_0'(y)}m_0(y)h_0(y)\,\mathrm{d}\mu(x)
\end{align*}
Thus
$\frac{1}{\#r^{-1}(x)}\sum_{r(y)=x}\cj{m_0'(y)}m_0(y)h_0(y)=h_0(x)$
$\mu$-a.e.
\par
Next we check that $|h_0|^2\leq\|A\|^2hh'$ $\mu$-a.e. By the
Schwarz inequality, we have for all $f,g\in L^\infty(X)$,
\begin{equation*}|\ip{\pi'(f)\varphi'}{A\pi(g)\varphi}|^2\leq\|A\|^2\|\pi'(f)\varphi'\|^2\|\pi(g)\varphi\|^2,\end{equation*}
which translates into
\begin{equation}\label{eqinterpf1}
|\int_X\cj
fgh_0\,\mathrm{d}\mu|^2\leq\|A\|^2\int_X|g|^2h'\,\mathrm{d}\mu\,\int_X|f|^2h\,\mathrm{d}\mu.
\end{equation}
If $\mu$ has some atoms then just take $f$ and $g$ to be the
characteristic function of that atoms and this proves the
inequality (\ref{eqinter3_0}) for such points. The part of $\mu$
that does not have atoms is measure theoretically isomorphic to
the unit interval with the Lebesgue measure. Then take $x$ to be a
Lebesgue differentiability point for $h_0,h$ and $h'$. Take
$f=g=\frac{1}{\mu(I)}\chi_I$ for some small interval centered at
$x$. Letting $I$ shrink to $x$ and using Lebesgue's
differentiability theorem, (\ref{eqinterpf1}) implies
(\ref{eqinter3_0}).
\par
For the converse, from (iii) to (i), let $h_0$ as in (iii), and
define for $n\geq0$ the sesquilinear form, $B_n$ on $H_n'\times
H_n$ (see Section \ref{Mul}): for $f,g\in L^\infty(X)$,
\begin{equation*}B_n(U'^{-n}\pi'(f)\varphi',U^{-n}\pi(g)\varphi):=\int_X\cj{f}gh_0\,\mathrm{d}\mu\end{equation*}
An application of the Schwarz inequality and (\ref{eqinter3_0}),
shows that
\begin{equation*}|B_n(\xi,\eta)|^2\leq c\|\xi\|^2\|\eta\|^2,\quad(\xi\in H_n',\eta\in H_n).\end{equation*}
The inclusion of $H_n$ in $H_{n+1}$ is given by
\begin{equation*}U^{-n}\pi(f)\varphi\mapsto U^{-n-1}\pi(f\circ r\,m_0)\varphi.\end{equation*}
The forms $B_n$ are compatible with these inclusion in the sense
that
\begin{multline*}B_{n+1}(U'^{-n-1}\pi'(f\circ rm_0')\varphi',U^{-n-1}\pi(g\circ rm_0)\varphi)\\
=\int_X\cj{f\circ rm_0'}g\circ rm_0h_0\,\mathrm{d}\mu=\int_X\cj fgh_0=B_n(U'^{-n}\pi'(f)\varphi',U^{-n}\pi(g)\varphi)\end{multline*}
(We used (\ref{eqinter3}) for the third equality.) Therefore the
system $(B_n)_n$ extends to a sesquilinear map $B$ on $H'\times H$
such that its restriction to $H_n'\times H_n$ is $B_n$, and $B$ is
bounded $(H=L^2(X_\infty,\hat\mu)$, $H'=L^2(X_\infty,\hat\mu')$.)
Then there exists a bounded operator $A:H\rightarrow H'$ such that
\begin{equation*}\ip{\xi}{A\eta}=B(\xi,\eta),\quad(\xi\in H,\eta\in H').\end{equation*}
We have to check that $A$ is intertwining. But
\begin{align*}
\bip{U'^{-n}\pi'(f)\varphi'}{AUU^{-n}\pi(g)\varphi}&=B(U'^{-n}\pi'(f)\varphi',U^{-n}\pi(g\circ r\,m_0)\varphi)\\
&=\int_X\cj{f}g\circ r\,m_0h_0\,\mathrm{d}\mu\\
&=B(U'^{-n-1}\pi'(f)\varphi',U^{-n-1}\pi(g\circ r\,m_0)\varphi)\\
&=\bip{U'^{-n-1}\pi'(f)\varphi'}{AU^{-n-1}\pi(g\circ r\,m_0)\varphi}\\
&=\bip{U'^{-n}\pi'(f)\varphi'}{U'AU^{-n}\pi(g)\varphi}.
\end{align*}
\begin{align*}
\bip{U'^{-n}\pi'(f)\varphi'}{A\pi(k)U^{-n}\pi(g)\varphi}&=B(U'^{-n}\pi'(f)\varphi',U^{-n}\pi(k\circ r^n\,g)\varphi)\\
&=\int_X\cj fk\circ r^n\,gh_0\,\mathrm{d}\mu\\
&=B(U'^{-n}\pi'(\cj{k\circ r^n}\,f)\varphi',U^{-n}\pi(g)\varphi)\\
&=\bip{U'^{-n}\pi'(\cj{k\circ r^n}\,f)\varphi'}{AU^{-n}\pi(g)\varphi}\\
&=\bip{\pi'(\cj k)U'^{-n}\pi'(f)\varphi'}{AU^{-n}\pi(g)\varphi}\\
&=\bip{U'^{-n}\pi'(f)\varphi'}{\pi'(k)AU^{-n}\pi(g)\varphi}.
\end{align*}
This shows that $A$ is intertwining.
\par
{}From (ii) to (iii), take $f$ as in (ii). Then define the operator
$A$ as in (\ref{eqinter4}). Using the previous correspondences we
have that $A$ is intertwining and there exists $h_0$ as in (iii),
satisfying (\ref{eqinter6}). We rewrite this in terms of $f$, and
we have for all $g\in L^\infty(X)$:
\begin{equation*}\int_XE_0^{\hat\mu'}(f)gh'\,\mathrm{d}\mu=\int_{X_\infty}fg\,\mathrm{d}\hat\mu'=\ip{\varphi'}{A\pi(g)\varphi}=\int_Xgh_0\,\mathrm{d}\mu\end{equation*}
Also
\begin{equation*}\int_{X_\infty}fg\,\mathrm{d}\hat\mu'=\int_{X_\infty^a}fg\Delta\,\mathrm{d}\hat\mu=\int_XE_0^{\hat\mu}(f\Delta)gh\,\mathrm{d}\mu.\end{equation*}
This proves (\ref{eqinter5}).
\end{proof}

\begin{corollary}\label{corinter2}
Let $(m_0,h)$ be a Perron-Ruelle-Frobenius pair with $m_0$
non-singular.
\begin{enumerate}
\item For each operator $A$ on $L^2(X_\infty,\hat\mu)$ which
commutes with $U$ and $\pi$, there exists a cocycle $f$, i.e., a
bounded measurable function $f:X_\infty\rightarrow\bc$ with
$f=f\circ\hat r$, $\hat\mu$-a.e., such that
\begin{equation}\label{eqinterc1}A=M_f,
\end{equation}
and, conversely each cocycle defines an operator in the commutant.
\item For each measurable harmonic function $h_0:X\rightarrow\bc$,
i.e., $R_{m_0}h_0=h_0$, with $|h_0|^2\leq ch^2$ for some $c\geq0$,
there exists a unique cocycle $f$ such that
\begin{equation}\label{eqinterc2}
h_0=E_0(f)h,
\end{equation}
and conversely, for each cocycle the function $h_0$ defined by
\textup{(\ref{eqinterc2})} is harmonic. \item The correspondence
$h_0\rightarrow f$ in \textup{(ii)} is given by
\begin{equation}\label{eqinterc3}
f=\lim_{n\rightarrow\infty}\frac{h_0}{h}\circ\theta_n
\end{equation}
where the limit is pointwise $\hat\mu$-a.e., and in
$L^p(X_\infty,\hat\mu)$ for all $1\leq p<\infty$.
\end{enumerate}
\end{corollary}
\begin{proof}
(i) and (ii) are direct consequences of Theorem \ref{thintertw}.
For (iii), we have that $f\in L^\infty(X_\infty,\hat\mu)\subset
L^p(X_\infty,\hat\mu)$. Using Proposition \ref{propmartins4_11},
we have that, since $f=f\circ\hat r$, if
$E_n(f)=\xi_n\circ\theta_n$, then
\begin{equation*}\xi_n=\xi_{n+1},\mbox{ for all }n\geq0.\end{equation*}
But from (\ref{eqinterc2}), we know that $\xi_0=\frac{h_0}{h}$, so
\begin{equation*}E_n(f)=\frac{h_0}{h}\circ\theta_n.\end{equation*}
(iii) follows now from Propositions \ref{propcond4},
\ref{propcond5} and \ref{propcond6}.
\end{proof}

\section{\label{Ite}\uppercase{Iterated function systems}}
In Section \ref{Mar} we constructed our extension systems using
martingales, and Doob's convergence theorem. We showed that our
family of martingale Hilbert spaces may be realized as
$L^2(X_\infty,\hat\mu)$, where both $X_\infty$, and the associated
measures $\hat\mu$ on $X_\infty$ are projective limits constructed
directly from the following given data. Our construction starts
with the following four: (1) a compact metric space $X$, (2) a
given mapping $r:X\rightarrow X$, (3) a strongly invariant measure
$\mu$ on $X$, and (4) a function $W$ on $X$ which prescribes
transition probabilities. From this, we construct our extension
systems.\par
   In this section, we take a closer look at the measure $\hat\mu$. We show that $\hat\mu$
is in fact an average over an indexed family of measures $P_x$,
$x$ in $X$. Now $P_x$ is constructed as a measure on a certain
space of paths. The subscript $x$ refers to the starting point of
the paths, and $P_x$ is defined on a sigma-algebra of subsets of
path-space. (The reader is referred to \cite{Jor06} for additional
details.)
\par These are paths of a random walk, and the random
walk is closely connected to the mathematics of the projective
limit construction in Section \ref{Mul}. But the individual measures $P_x$
carry more information than the averaged version $\mu$ from
Section \ref{Mul}. As we show below, the construction of solutions to the
canonical scaling identities in wavelet theory, and in dynamics,
depend on the path space measures $P_x$. Our solutions will be
infinite products, and the pointwise convergence of these infinite
products depends directly on the analytic properties of the
$P_x$'s.
\par
Let $X$ be a metric space and $r:X\rightarrow X$ an $N$ to $1$
map. Denote by $\tau_k:X\rightarrow X$, $k\in\{1,...,N\}$, the
branches of $r$, i.e., $r(\tau_k(x))=x$ for $x\in X$, the sets
$\tau_k(X)$ are disjoint and they cover $X$.
\par
Let $\mu$ be a measure on $X$ with the property
\begin{equation}\label{eq3_1}
\mu=\frac{1}{N}\sum_{k=1}^N\mu\circ\tau_k^{-1}.
\end{equation}
This can be rewritten as
\begin{equation}\label{eq3_2}
\int_{X}f(x)\,\mathrm{d}\mu(x)=\frac{1}{N}\sum_{k=1}^N\int_Xf(\tau_k(x))\,\mathrm{d}\mu(x),
\end{equation}
which is equivalent also to the strong invariance property. \par

Let $W, h\geq0$ be two functions on $X$ such that
\begin{equation}\label{eq3_3}
\sum_{k=1}^NW(\tau_k(x))h(\tau_k(x))=h(x),\quad(x\in X).
\end{equation}
Denote by
\begin{equation*}\Omega:=\Omega_N:=\prod_{\mathbb{N}}\{1,...,N\}.\end{equation*} Also we denote by
\begin{equation*}W^{(n)}(x):=W(x)W(r(x))...W(r^{n-1}(x)),\quad(x\in X).\end{equation*}
\begin{proposition}\label{prop3_1}
For every $x\in X$ there exists a positive Radon measure $P_x$ on
$\Omega$ such that, if $f$ is a bounded measurable function on
$\Omega$ which depends only on the first $n$ coordinates
$\omega_1,...,\omega_n$, then
\begin{multline}\label{eq3_4}\int_\Omega
f(\omega)\,\mathrm{d}P_x(\omega)\\
=\sum_{\omega_1,...,\omega_n}W^{(n)}(\tau_{\omega_n}\tau_{\omega_{n-1}}...\tau_{\omega_1}(x))h(\tau_{\omega_n}\tau_{\omega_{n-1}}...\tau_{\omega_1}(x))f(\omega_1,...,\omega_n).
\end{multline}
\end{proposition}
\begin{proof}
We check that $P_x$ is well defined on functions which depend only
on a finite number of coordinates. For this, take $f$ measurable
and bounded, depending only on $\omega_1,...,\omega_n$; and
consider it as function which depends on the first $n+1$
coordinates. We have to check that the two formulas given by
(\ref{eq3_4}) yield the same result.
\par
{\it Consistency}: As a function of the first $n+1$ coordinates,
we have
\begin{align*}
\int_Xf(\omega)\,\mathrm{d}P_x(\omega)&=\sum_{\!\omega_1,...,\omega_{n+1}\!}W^{(n+1)}(\tau_{\omega_{n+1}}...\tau_{\omega_1}(x))h(\tau_{\omega_{n+1}}...\tau_{\omega_1}(x))f(\omega_1,...,\omega_{n+1})\\
&=\sum_{\omega_1,...,\omega_n}f(\omega_1,...,\omega_n)W^{(n)}(\tau_{\omega_{n}}...\tau_{\omega_1}(x))\\ &\qquad{}\cdot\sum_{\omega_{n+1}} W(\tau_{\omega_{n+1}}...\tau_{\omega_1}(x))h(\tau_{\omega_{n+1}}...\tau_{\omega_1}(x))\\
&=\sum_{\omega_1,...,\omega_n}W^{(n)}(\tau_{\omega_{n}}...\tau_{\omega_1}(x))f(\omega_1,...,\omega_n)h(\tau_{\omega_{n}}...\tau_{\omega_1}(x)),
\end{align*}
so $P_x$ is well defined. Using the Stone-Weierstrass and Riesz
theorems, we obtain the desired measure.
\end{proof}
Consider now the space $X\times\Omega$. On this space we have the
shift $S$:
\begin{equation}\label{eq3_5}
S(x,\omega_1...\omega_n...)=(r(x),\omega_x\omega_1...\omega_n...),\quad(x\in
X,(\omega_1...\omega_n...)\in\Omega),
\end{equation}
where $\omega_x$ is defined by $x\in\tau_{\omega_x}(X)$. The
inverse of the shift is given by the formula:
\begin{equation}\label{eq3_6}
S^{-1}(x,\omega_1...\omega_n...)=(\tau_{\omega_1}(x),\omega_2...\omega_n...),\quad(x\in
X,(\omega_1...\omega_n...)\in\Omega).
\end{equation}
\begin{proposition}\label{prop3_2}
Define the map $\Psi:X_\infty\rightarrow X\times\Omega$ by
\begin{equation*}\Psi(x_0,x_1,...)=(x_0,\omega_1,\omega_2,...),\mbox{ where }x_{n}=\tau_{\omega_{n}}(x_{n-1}),\quad(n\geq 1).\end{equation*}
Then $\Psi$ is a measurable bijection with inverse
\begin{equation*}\Psi^{-1}(x,\omega_1,\omega_2,...)=(x,\tau_{\omega_1}(x),\tau_{\omega_2}\tau_{\omega_1}(x),...).\end{equation*}
\begin{equation}\label{eq3_2_0}
\Psi\circ\hat r\circ\Psi^{-1}=S. \end{equation} Also
\begin{equation}\label{eq3_2_1}
\int_{X_\infty}f\,\mathrm{d}\hat\mu=\int_X\int_\Omega
f\circ\Psi^{-1}(x,\omega)\,\mathrm{d}P_x(\omega)\,\mathrm{d}\mu(x),\quad(f\in
L^1(X_\infty,\hat\mu)).
\end{equation}
\end{proposition}
\begin{proof}
We know that $r(x_{n})=x_{n-1}$ therefore
$x_{n}=\tau_{\omega_n}(x_{n-1})$ for some
$\omega_n\in\{1,...,N\}$. This correspondence defines $\Psi$ and
it is clear that the map is 1-1 and onto and the inverse has the
given formula. A computation proves (\ref{eq3_2_0}).
\par
To check (\ref{eq3_2_1}), it is enough to verify the conditions of
Theorem \ref{thpr_1}. Take $\xi\in L^\infty(X)$, then
$\xi\circ\theta_n\circ\Psi^{-1}$ depends only on $x$ and
$\omega_1,...,\omega_n$ so
\begin{align*}&\int_X\int_\Omega
f\circ\theta_n\circ\Psi^{-1}(x,\omega)\,\mathrm{d}P_x(\omega)\,\mathrm{d}\mu(x)\\
&\qquad=\int_X\sum_{\omega_1,...,\omega_n}W^{(n)}(\tau_{\omega_n}\tau_{\omega_{n-1}}...\tau_{\omega_1}(x))\\
&\qquad\qquad\quad{}\cdot h(\tau_{\omega_n}\tau_{\omega_{n-1}}...\tau_{\omega_1}(x))(f\theta_n\Psi^{-1})(x,\omega)\,\mathrm{d}\mu(x)\\
&\qquad=\int_XR^n(fh)(x)\,\mathrm{d}\mu(x)=\int_{X_\infty}f\circ\theta_n\,\mathrm{d}\hat\mu.\end{align*}
This proves (\ref{eq3_2_1}).
\end{proof}

\begin{acknowledgements}
We gratefully acknowledge discussions
with the members of the group FRG/NSF; especially enlightening
suggestions from Professors L. Baggett, D. Larson, and G.
Olafsson.

In addition, this work was supported at the University of Iowa by
a grant from the National Science Foundation (NSF-USA) under a
Focused Research Program, DMS-0139473 (FRG).
\end{acknowledgements}

\enlargethispage{0.5\baselineskip}
\end{document}